\documentclass[final,12pt]{elsarticle}

\usepackage[letter paper, margin=1 in]{geometry}
\usepackage[numbers]{natbib}

\usepackage[utf8]{inputenc} 
\usepackage[T1]{fontenc}    
\usepackage[hidelinks,colorlinks=true,linkcolor=blue,citecolor=blue]{hyperref}
\usepackage{url}            
\usepackage{booktabs}       
\usepackage{dsfont}         
\usepackage{algorithm}
\usepackage{algorithmic}
\usepackage[algo2e,ruled,vlined,linesnumbered]{algorithm2e}

\usepackage{amsfonts,amsmath,amssymb,amsthm}      
\usepackage{bbm}
\usepackage{nicefrac}       
\usepackage{microtype}      
\usepackage[capitalize,noabbrev]{cleveref}
\usepackage{graphicx}
\usepackage{wrapfig}
\usepackage{subcaption}
\usepackage{xcolor}
\usepackage{comment}
\usepackage{empheq}
\usepackage[most]{tcolorbox}

\newcommand{\mbb}{\mathbb}
\newcommand{\mbf}{\mathbf}
\newcommand{\mcl}{\mathcal}
\newcommand{\bs}{\boldsymbol}

\newcommand{\x}{\mathbf{x}}
\newcommand{\X}{\mathcal{X}}
\newcommand{\D}{\bs{\mathcal{D}}}

\newcommand{\qpots}{$q\texttt{POTS}$}
\DeclareMathOperator*{\argmax}{argmax}

\newcommand{\mg}{\textcolor{black}}
\theoremstyle{plain}
\newtheorem{theorem}{Theorem}[section]
\newtheorem{proposition}[theorem]{Proposition}

\theoremstyle{definition}

\theoremstyle{remark}

\newtcbox{\mymath}[1][]{%
    nobeforeafter, math upper, tcbox raise base,
    enhanced, colframe=gray!30!black,
    colback=gray!30, boxrule=1pt,
    #1}
\tcbset{highlight math style={boxsep=5mm,colback=gray!30!red!30!white}}

\title{Multiobjective Aerodynamic Design Optimization of the NASA
Common Research Model}

\author[psu]{Kade Carlson}
\author[psu]{Ashwin Renganathan}

\address[psu]{The Pennsylvania State University, \\
556 White Course Drive, University Park, PA 16802, USA}

\begin{document}

\begin{frontmatter}


\begin{abstract}
Aircraft aerodynamic design optimization must account for the varying operating conditions along the cruise segment as opposed to designing at one fixed operating condition, to arrive at more realistic designs. Conventional approaches address this by performing a ``multi-point'' optimization that assumes a weighted average of the objectives at a set of sub-segments along the cruise segment. We argue that since such multi-point approaches are, inevitably, biased by the specification of the weights, they can lead to sub-optimal designs. Instead, we propose to optimize the aircraft design at multiple sub-segments simultaneously -- that is, via multiobjective optimization that leads to a set of Pareto optimal solutions. However, existing work in multiobjective optimization suffers from (i) lack of sample efficiency (that is, keeping the number of function evaluations to convergence minimal), (ii) scalability  \mg{in the absence of derivative information}, and (iii) the ability to generate a batch of iterates for synchronous parallel evaluations. To overcome these limitations, we \mg{apply} a novel multiobjective Bayesian optimization methodology \mg{for aerodynamic design optimization} that demonstrates improved sample efficiency  and accuracy compared to the state of the art. Inspired by Thompson sampling, our approach leverages Gaussian process surrogates and Bayesian decision theory to generate a sequence of iterates according to the probability that they are Pareto optimal. Our approach, named batch Pareto optimal Thompson sampling (\qpots)\footnote{Here, $q$ stands for selecting a batch of $q$ iterates at every step.}, demonstrates superior empirical performance on a variety of synthetic experiments as well as a $24$ dimensional two-objective aerodynamic design optimization of the NASA common research model. We also provide open-source software of our methodology \mg{and experiments}.
\end{abstract}




\end{frontmatter}



\section{Introduction}
\label{sec:intro}
\begin{figure}[htb!]
    \centering
    \includegraphics[width=1\linewidth]{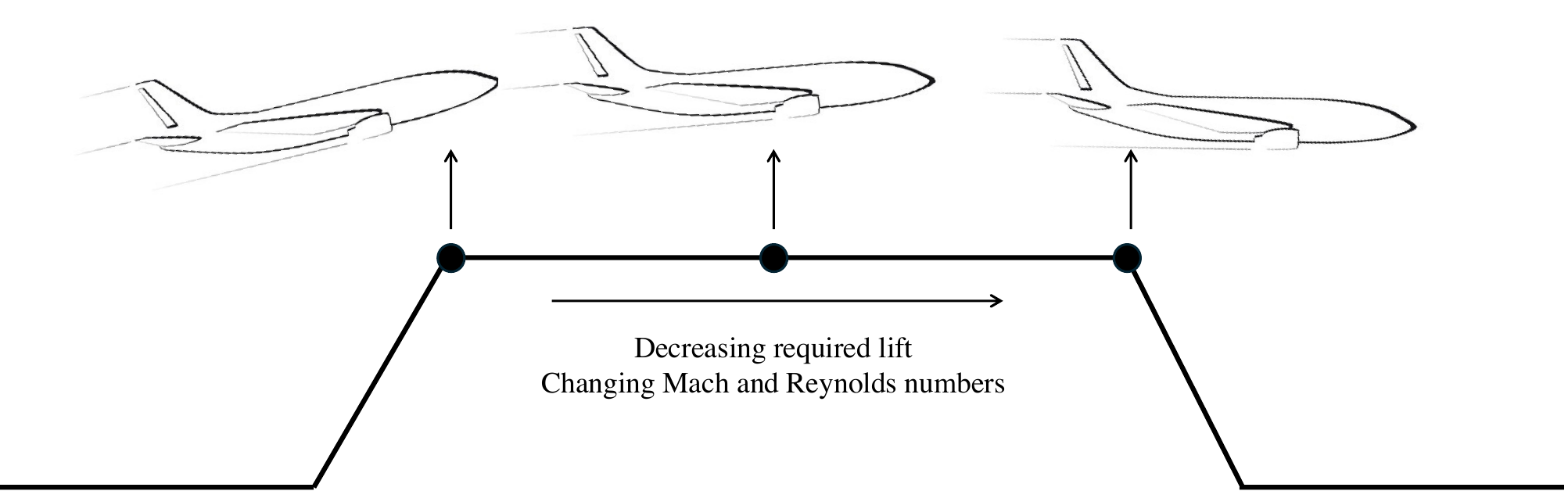}    
    \caption{A simplified mission profile for a commercial aircraft that emphasizes the cruise segment. The change in aircraft weight during the cruise segment induces change in required lift and thrust, and hence the operating angle of attack, Mach and Reynolds numbers. We anticipate that designs optimized for any one point along the cruise segment will likely be sub-optimal for other points. Therefore, we argue that ``multiobjective'' optimization, that simultaneously considers multiple points along the cruise segment is necessary.}
    \label{fig:cruise_pts}
\end{figure}

Aerodynamic design optimization refers to optimizing the wing and control surface cross-sectional shapes, typically to minimize drag or maximize the lift-to-drag ratio (aka, the aerodynamic efficiency). As a result, it forms a routine part of the detailed design phase of aircraft design when the major design decisions are locked in. Conventionally, the aircraft wing is optimized for the cruise condition, which is reasonable given that it is the most dominant segment of a commercial aircraft's mission. In this regard, it is quite common to consider a constant operating condition for the cruise segment. In principle, this is not a valid assumption because the aircraft weight is constantly decreasing along this segment (due to fuel consumption), resulting in changing required thrust, which in turn changes the operating angle of attack, Mach and Reynolds numbers; we provide a graphical illustration in \Cref{fig:cruise_pts}.  

\begin{figure}[h!]
    \centering
    \begin{subfigure}{.5\textwidth}
        \includegraphics[width=1\linewidth]{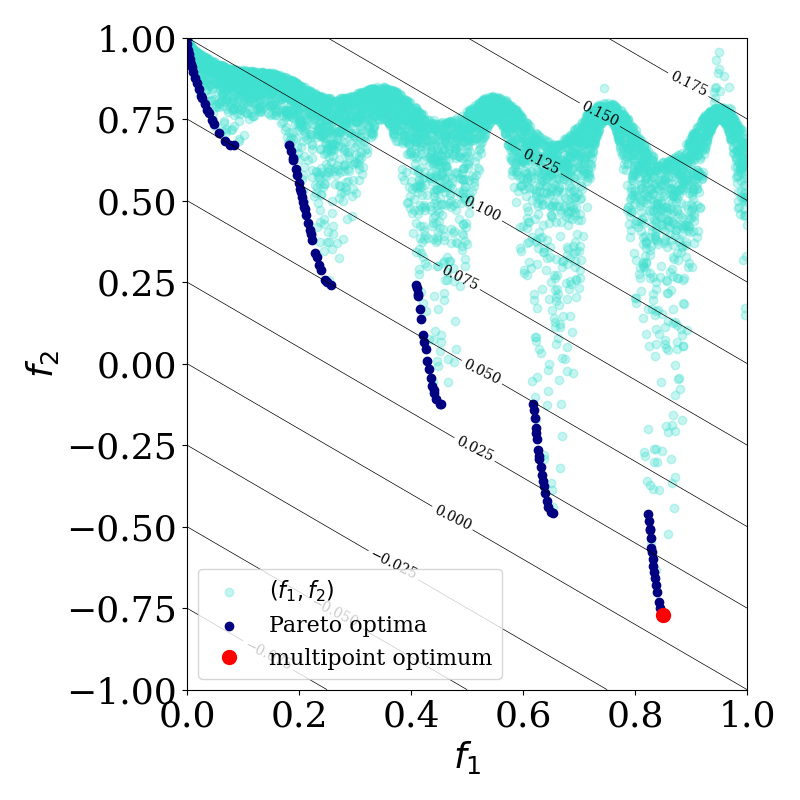}
    \end{subfigure}%
    \begin{subfigure}{.5\textwidth}
        \includegraphics[width=1\linewidth]{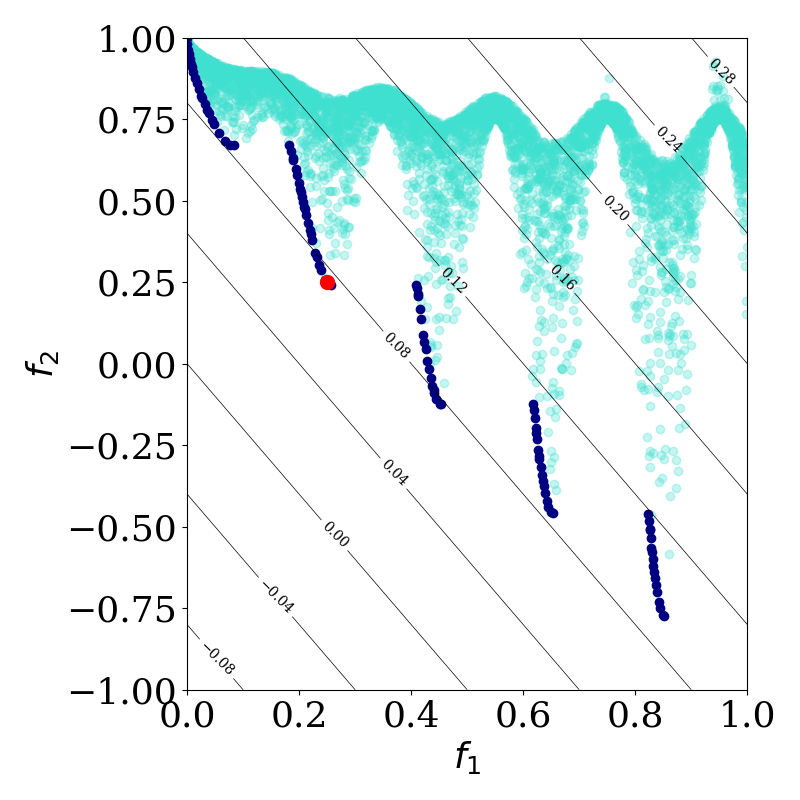}
    \end{subfigure}
    \caption{Illustration of the differences between multi-point versus multiobjective (Pareto) optima on the two-objective ZDT3 test function. The turqouise circles are the objective pairs and the blue circles are the Pareto optimal objectives. The lines with inline labels represent the multi-point (scalarized) combination of the objectives: $w_1 f_1 + w_2 f_2$, for two different choices of weights $w_1,w_2$ on the left and right. The red circle is the multi-point optimum which is essentially the point  where the lines intersect the Pareto frontier with the smallest \emph{scalarized} objective. As can be seen, the choice of weights can drastically change the multi-point optimum. On the other hand, multiobjective optimization provides an ensemble of optimal designs that a practitioner can choose from.}
    \label{fig:zdt3}
\end{figure}

In principle, one must consider a range of angles of attack, Mach, and Reynolds numbers in optimizing the wing shape. The aerodynamic design optimization discussion group (ADODG)~\footnote{\url{https://sites.google.com/view/mcgill-computational-aerogroup/adodg}}, accordingly, propounds $3$ and $9$ point ``multi-point`` design optimization, which considers a weighted average over a set of operating conditions from the cruise segment. Those problems, however, involve a predefined weight for each sub-segment, which would need some domain knowledge to be specified appropriately. In this work, we argue that each sub-segment for the cruise segment must be simultaneously optimized, as opposed to scalarizing them using a bunch of weights. This is for the following reasons. First, we anticipate that the designs at individual sub-segments are conflicting with each other -- that is, the optimal design at one segment can be sub-optimal at another. Second, the choice of weights, inevitably, affects the final design and hence must be carefully chosen, which might not be practical. Finally, simultaneously optimizing all objectives recovers an ensemble of solutions, of which the scalarized ``multi-point'' solution is one. 
\mg{Specifically, minimizing a weighted sum with strictly positive weights always yields a Pareto–optimal solution
(hence a scalarized solution is on the Pareto front). However, the converse requires convexity of
the attainable objective set: in nonconvex problems linear scalarizations cannot recover Pareto
points in concave regions of the front, which motivates the use of 
direct multiobjective methods like ours.}
 We illustrate this on a synthetic test function, the ZDT3 test function, in \Cref{fig:zdt3}, and provide a concrete theoretical justification in \Cref{sec:multip_v_multiob}. 

\begin{figure}[htb!]
    \centering
    \includegraphics[width=0.75\linewidth]{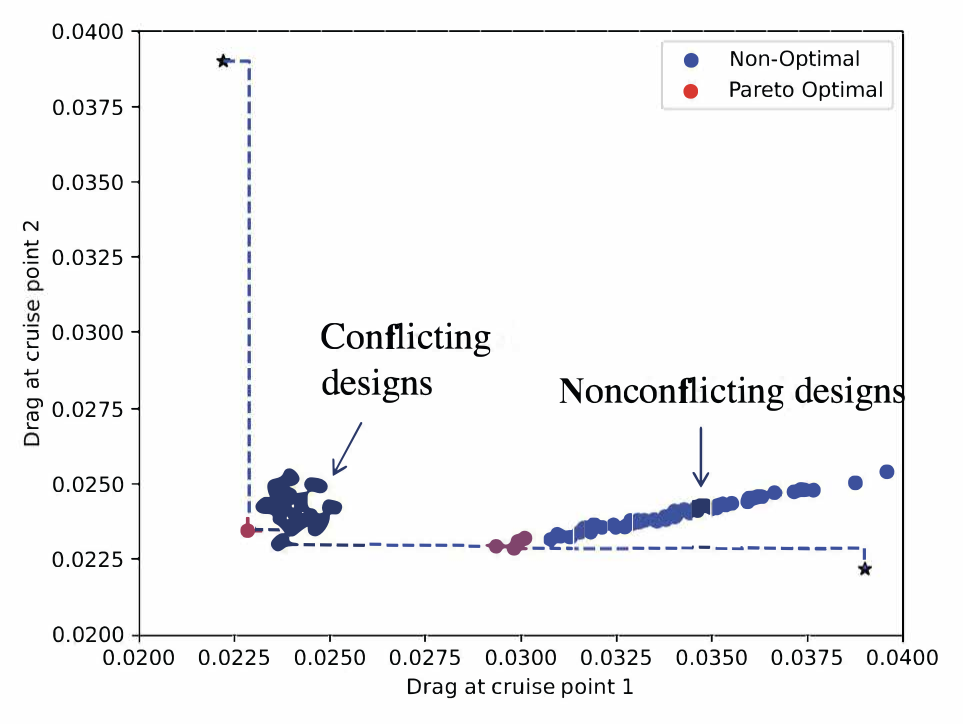}
    \caption{Plot shows aircraft drag for the NASA CRM evaluated at two different points along cruise segment for a variety of wing shapes. Notice that certain designs are nonconflicting, while others are not. Multiobjective optimization is necessary to properly account for design trade-offs and identify a set of optimal solutions for a practitioner to choose from.}
    \label{fig:crm_pareto}
\end{figure}

Such a simultaneous optimization of a set of conflicting objectives is called multiobjective optimization, which forms the focus of this work -- we provide a formal introduction in \Cref{sec:problem_statement}. This is crucial for aircraft design because the design space can include regions of both \mg{conflicting} and \mg{nonconflicting} designs; this is shown in \Cref{fig:crm_pareto}, where we consider the aircraft drag at two different cruise points, which are to be minimized. When designs are \mg{nonconflicting} (as in the bottom right of \Cref{fig:crm_pareto}), then a multipoint optimum, regardless of the choice of weights, will find the same optimal design. However, when the designs are \mg{conflicting}, a ``Pareto optimal'' set of designs must be recovered (red circles in \Cref{fig:crm_pareto}). Therefore, efficient methods to simultaneously optimize objectives at multiple operating conditions are crucial -- this work makes contributions in that direction.

Multiobjective optimization has been widely adopted for aerodynamic design optimization in the past. Evolutionary algorithms, particularly genetic algorithms, are arguably the most popular. For instance, \citet{vicini_inverse_1997} present a numerical procedure for transonic airfoil design using genetic algorithms with single-point, multi-point, and multiobjective optimization formulations. The study highlights the advantages of multiobjective optimization over single-point and multi-point approaches in optimizing drag rise characteristics of transonic airfoils. Similarly, \citet{coello_coello_moses:_1999} introduce two new multiobjective optimization techniques based on genetic algorithms, demonstrating their ability to generate better trade-offs compared to other approaches. Advanced techniques such as scatter search and adaptive range multiobjective genetic algorithms have been applied to complex aerodynamic design problems. \citet{molina_sspmo:_2007} describe a metaheuristic procedure based on scatter search for approximating the efficient frontier of nonlinear multiobjective optimization problems, showing its viability as an alternative for multiobjective optimization. \citet{chiba_multidisciplinary_2007} report a large-scale application of evolutionary multiobjective optimization to the multidisciplinary design optimization of a transonic regional-jet aircraft wing, employing high-fidelity evaluation models and an adaptive range multiobjective genetic algorithm.

Another angle is hybrid approaches combining evolutionary algorithms with surrogate models that can offer a sample-efficient (that is, minimal number of evaluations of the objectives) approach to multiobjective optimization. \citet{leifsson_multiobjective_2015} present a computationally efficient procedure for multiobjective design optimization using variable-fidelity models and response surface surrogates. This approach combines an evolutionary algorithm with Kriging interpolation and space-mapping correction, allowing for efficient exploration-exploitation of trade-offs between conflicting objectives in airfoil design. \citet{kim_axial-flow_2011} also demonstrates the effectiveness of a hybrid multiobjective evolutionary algorithm coupled with a response surface approximation surrogate model for axial-flow ventilation fan design, achieving significant improvements in efficiency and pressure rise. \citet{keane_statistical_2006} presents an approach that enables the parallel evolution of multiobjective Pareto sets, providing a statistically coherent means of solving expensive multiobjective design problems. \mg{Other recent efforts include \cite{zuhal2019comparative} and \cite{lin2023parallel}}.

Despite the advances made in the past, multiobjective aerodynamic design optimization still requires further attention. First, existing algorithms, such as evolutionary approaches, can incur a large number of objective evaluations, which can be prohibitive when using expensive high-fidelity aerodynamic models. Second, the efficient algorithms that exist today require derivative information, which can be difficult or impossible to obtain in several practical scenarios; see \cite{plaban2025derivative} for a detailed discussion. Finally, in existing work \mg{sampling a batch of iterates at every step, for asynchronous parallel evaluations on high-performance computing environments, is either not possible or incurs a nontrivial additional cost}. This work aims to address these limitations. Specifically, we take a surrogate-based approach, where we fit global surrogate models for our objectives. Then, via Bayesian decision theory, we adaptively enrich the surrogate model until a suitable convergence criterion is satisfied -- this, in turn, results in a desirable approximation to the solution to the multiobjective optimization problem. Overall, our approach falls under the Bayesian optimization (BO)~\cite{frazier2018tutorial,shahriari2015taking} family of methods, and we provide more detail in \Cref{sec:method}. BO has been successfully demonstrated in the past to be competitive with the classical adjoint-based methods for aerodynamic design optimization~\cite{renganathan2021enhanced, plaban2025derivative}.

\mg{In this work, we argue that multiobjective aerodynamic design optimization is more effective than solving multiple multi-point optimization problems. This is because, a single solution of a multiobjective optimization problem results in a set of designs that covers every possible multi-point design. The converse, however, is not true: multi-point optimization, regardless of the scalarizing weights, is not guaranteed to find a Pareto optimal point. We substantiate this argument theoretically in \Cref{sec:multip_v_multiob}. Therefore, naturally, direct multiobjective optimization methods are necessary. In this regard, to address the existing limitations in solving multiobjective aerodynamic design optimization problems in the absence of adjoint information, we propose a novel multiobjective Bayesian optimization method and demonstrate it on the aerodynamic design optimization of the NASA common research model (CRM).
\newline
\noindent {\bf Remark 1.} As a first step, our study is scoped down to unconstrained multiobjective aerodynamic design optimization of the NASA CRM under inviscid flow conditions; we briefly discuss constrained extensions in \Cref{sec:constraints}. However, we rigorously benchmark the proposed unconstrained optimization method against the state of the art to establish the merits of the method.
\newline
\noindent {\bf Remark 2.} In the current work, we explore multiobjective optimization of modest input ($2 - 24$) and output dimensionalities ($2$); high-dimensional scaling is reserved for future work. Demonstration on fully viscous constrained multiobjective aerodynamic design is reserved for future work.
}

The rest of the paper is organized as follows. In \Cref{sec:problem_statement} we formally introduce our multiobjective optimization problem of interest. In \Cref{sec:method} we present our overall methodology, and in \Cref{sec:experiments} we provide details on our experiments and discuss results. Finally, we conclude in \Cref{sec:conclusion} with an outlook for future work.

\section{Problem statement and preliminaries}
\label{sec:problem_statement}
Mathematically, we consider the simultaneous unconstrained optimization of $K$ objectives 
\begin{equation} 
    \begin{split}
    \max_{\x \in \X}~&\{f_1(\x), \ldots, f_K(\x) \}     
    \end{split}
    \label{eqn:main_problem}
    \end{equation}
where $\x \in \X \subset \mbb{R}^d$ is the design variable, $f_k: \X \rightarrow \mbb{R},~\forall k=1,\ldots,K$, are expensive zeroth-order oracles (i.e., no available derivative information). $\X$ is the domain of $f$. We are interested in identifying a \emph{Pareto} optimal set of solutions, which ``Pareto dominates'' all other points. Let $\bs{f(\x)}=[f_1(\x),\ldots,f_K(\x)]^\top$; when a solution $\bs{f(\x)}$ Pareto dominates another solution $\bs{f(\x')}$, then $f_k(\x) \geq f_k(\x'),~\forall k=1,\ldots,K$ and $ \exists k \in [K]$ such that $f_k(\x) > f_k(\x')$. We write Pareto dominance as $\bs{f(\x)} \succ \bs{f(\x')}$. The set $\mcl{Y}^* = \{\bs{f(\x)} ~:~ \nexists \x' \in \X: \bs{f(\x') \succ \bs{f(\x)}} \}$ is called the \emph{Pareto frontier} and the set $\X^* = \{ \x \in \X ~:~ \bs{f(\x)} \in \mcl{Y}^*\}$ is called the \emph{Pareto set}.

 We focus on sample-efficient methods to address our problem. Evolutionary algorithms (EAs), such as the nondominated sorting genetic algorithm (NSGA)\cite{deb2000fast, deb2002fast}, are widely used for multiobjective optimization (see \cite{zitzler2000comparison} for a review). However, EA-based methods typically require a large number of evaluations, which can be impractical in settings involving expensive function evaluations. To overcome this limitation, a common strategy is to employ Bayesian optimization (BO) with Gaussian process (GP) surrogate models for each objective. In this approach, a GP is trained for each objective using an initial dataset $\D^i_n = {(\x_j, y_j),j=1,\ldots,n}$, where $y_j = f_i(\x_j)$. An \emph{acquisition} function is then constructed to measure the utility of a candidate point $\x$, and a subproblem is solved to identify the point(s) that maximize this acquisition function. The surrogate models are subsequently updated with observations at the new points, and this process is iterated until a chosen convergence condition is satisfied. We formally introduce this methodology in \Cref{sec:bo}. BO is especially popular in the single-objective case ($K=1$); for example, the expected improvement (EI)\cite{jones1998efficient} and probability of improvement (PI)\cite{mockus1978application,jones2001taxonomy} quantify the likelihood that a candidate point improves upon the best observed value. The GP upper confidence bound (UCB)\cite{srinivas2009gaussian} offers a conservative estimate of the maximum by leveraging the GP’s mean and uncertainty estimates. Other notable acquisition strategies for $K=1$ include entropy-based methods\cite{wang2017max,hernandez2014predictive}, the knowledge gradient (KG)\cite{frazier2008knowledge,wu2016parallel,wu2019practical}, and stepwise uncertainty reduction (SUR)\cite{picheny2014stepwise,chevalier2014fast}. Multiobjective Bayesian optimization (MOBO) generalizes BO to handle multiple objectives; we review relevant prior work in the sections that follow.

\subsection{Single objective Bayesian optimization}
\label{sec:bo}

We place a Gaussian process (GP) prior on the oracle, $f(\x) \sim \mcl{GP}(0, k(\x, \cdot))$, where $k(\cdot, \cdot): \X \times \X \rightarrow \mbb{R}+$ denotes the covariance function (also known as the \emph{kernel}). Observations from the oracle are modeled as $y_i = f(\x_i) + \epsilon_i$ for $i=1,\ldots,n$, where each $\epsilon_i$ is assumed to be zero-mean Gaussian noise with unknown variance $\tau^2$. Using the observed dataset $\D_n = \{\x_i, y_i, \tau^2\}_{i=1}^n$, we fit a posterior GP—parameterized by hyperparameters $\bs{\Omega}$—to obtain the following conditional posterior distribution~\cite{rasmussen:williams:2006}:
\begin{equation}
\begin{split}
Y(\x) \mid \D_n, \bs{\Omega} &\sim \mcl{GP}(\mu_n(\x), \sigma^2_n(\x)), \\
\mu_n(\x) &= \mbf{k}_n^\top [\mbf{K}_n + \tau^2 \mbf{I}]^{-1} \mbf{y}_n, \\
\sigma^2_n(\x) &= k(\x, \x) - \mbf{k}_n^\top [\mbf{K}_n + \tau^2 \mbf{I}]^{-1} \mbf{k}_n,
\end{split}
\label{e:GP}
\end{equation}
where $\mbf{k}_n \equiv k(\x, X_n)$ is the vector of covariances between the point $\x$ and all observation sites in $\D_n$, $\mbf{K}_n \equiv k(X_n, X_n)$ is the covariance matrix over the observed points, $\mbf{I}$ is the identity matrix, and $\mbf{y}_n$ is the vector of corresponding outputs. The observation sites are denoted by $X_n = [\x_1,\ldots,\x_n]^\top \in \mbb{R}^{n \times d}$. Bayesian optimization  then proceeds by making sequential decisions using a utility function defined over the posterior GP: $u(\x) = u(Y(\x)\mid \D_n)$. Each decision results from solving an "inner" optimization subproblem that maximizes an acquisition function, typically of the form $\alpha(\x) = \mbb{E}_{Y \mid \D_n}[u(Y(\x))]$.

\paragraph{Thompson sampling.}
An alternative to constructing acquisition functions is Thompson sampling (TS)\cite{thompson1933likelihood}, which selects actions based on the probability that they are optimal in a specified sense\cite{russo2014learning}. In the context of single-objective optimization, TS selects points according to the distribution $p_{\x^*}(\x)$, where $\x^*$ denotes the maximizer of the objective function. Within BO, this distribution can be written as:
\begin{equation}
\begin{split}
p_{\x^*}(\x) =& \int p_{\x^*}(\x \mid Y) p(Y \mid \D_n) dY 
= \int \delta\left(\x - \arg\max_{\x \in \X} Y(\x)\right) p(Y \mid \D_n) dY,
\end{split}
\label{eqn:ts}
\end{equation}
where $\delta$ denotes the Dirac delta function. In other words, we construct a distribution by concentrating the mass on every maximizer of the GP $Y$, and integrating it out by the posterior GP distribution.
Sampling from $p_{\x^*}(\x)$ thus corresponds to drawing a sample path $Y(\cdot, \omega)$—with $\omega$ representing the random realization—and finding its maximizer. Consequently, in Thompson sampling for BO, the acquisition function is simply the posterior sample path itself, $\alpha \equiv Y$. Our proposed method builds on TS, which we now present.

\section{Methodology}
\label{sec:method}

\subsection{Multiobjective Bayesian optimization}
\label{ss:mobo}
In contrast to single objective BO, here we first fit $K$ independent posterior GP models: $\mcl{M}_1, \ldots, \mcl{M}_K$ for the $K$ objectives, with observations $\D^{1:K}_n \triangleq \{ \D_n^1,\ldots, \D_n^K \}$.
In the multiobjective BO, acquisition functions are typically defined in terms of a hypervolume (HV) indicator. HV of a Pareto frontier $\mcl{Y}^*$ is, loosely speaking, the volume bounded by the Pareto frontier and a reference point $\mbf{r} \in \mbb{R}^K$; $\mbf{r}$ is typically chosen as the independent infimums of each objective (in the maximization context).
Expected hypervolume improvement (EHVI)~\cite{couckuytFastCalculationMultiobjective2014, emmerichComputationExpectedImprovement2008} computes the expectation (with respect to the posterior GP) of the improvement in HV due to a candidate point. The HV computation is nontrivial and scales poorly with number of objectives. Our approach, however, does not involve the computation of HV but is a direct extension of \cref{eqn:ts} to the multiobjective setting.

\subsection{$q$\texttt{POTS:} Batch Pareto optimal Thompson sampling for MOBO}
\label{s:qpots}

Inspired by Thompson sampling (TS), our method—\emph{batch Pareto optimal Thompson sampling} (\qpots)~\cite{renganathan2025q}—selects points based on the probability that they are Pareto optimal. Specifically, we sample $\x$ according to the probability
\begin{equation}
    \begin{split}
    p_{\X^*}(\x) = \int \delta \left(\x - \arg\max_{\x \in \X} \{Y_1(\x), \ldots, Y_K(\x)\} \right) \, p(Y_1 \mid \D^1_n) \cdots p(Y_K\mid\D^K_n) \, dY_1 \cdots dY_K,
    \end{split}
\end{equation}
where the integration is over the joint support of the posterior GPs for all objectives. In practice, this corresponds to sampling posterior paths $Y_i(\cdot, \omega)$ for each objective $i = 1, \ldots, K$ and selecting points from their Pareto set:
\[
\x_{n+1} \in \arg\max_{\x \in \X} \{Y_1(\x, \omega), \ldots, Y_K(\x, \omega)\}.
\]
To draw sample paths from the posterior, we use the reparameterization trick (omitting indexation over objectives for brevity):  
\[
Y(\cdot, \omega) = \mu_n(\cdot) + \mbf{\Sigma}_n^{1/2}(\cdot) \cdot Z(\omega),
\]  
where $Z$ is a standard multivariate normal random variable. We then solve the following multiobjective optimization problem:
\begin{equation}
    X^* = \arg\max_{\x \in \X} \{Y_1(\x, \omega), \ldots, Y_K(\x, \omega)\},
    \label{eqn:moo_gp}
\end{equation}
which can be tackled with any multiobjective solver. We use evolutionary methods (e.g., NSGA-II~\cite{deb2002fast}) to combine the sample efficiency of MOBO with the global search capabilities of evolutionary algorithms.

From the predicted Pareto set $X^*$, we then select a batch of $q$ points $\x_{n+1:q} \subset X^*$. Each point in this batch is Pareto optimal with respect to the sampled posterior paths. Although any point from $X^*$ is valid in principle, we promote diversity by selecting acquisitions according to a \emph{maximin distance} criterion~\cite{johnson1990minimax,sun2019synthesizing}. This encourages a natural trade-off between exploration and exploitation.

Let $\gamma(\cdot, \cdot): \mbb{R}^d \times \mbb{R}^d \to \mbb{R}_+$ denote the Euclidean distance. The batch of $q$ new points is then selected via the following sequence of maximin subproblems:
\begin{equation}
    \begin{split}
        \x_{n+1} &\in \arg\max_{{\x}^* \in X^*} \min_{\x_i \in X_n} \gamma(\x^*, \x_i), \\
        \x_{n+2} &\in \arg\max_{{\x}^* \in X^*} \min_{\x_i \in X_n \cup \{\x_{n+1}\}} \gamma(\x^*, \x_i), \\
        \vdots \\
        \x_{n+q} &\in \arg\max_{{\x}^* \in X^*} \min_{\x_i \in X_n \cup \{\x_{n+1}, \ldots, \x_{n+q-1}\}} \gamma(\x^*, \x_i).
    \end{split}
    \label{eqn:maximin}
\end{equation}

While \Cref{eqn:maximin} may appear computationally demanding, it is efficient in practice. This is because the Pareto set $X^*$ returned by evolutionary solvers has finite cardinality. As a result, we only need to compute an $N^* \times n$ pairwise distance matrix (where $N^* = |X^*|$) and sort the resulting vector of column-wise minima. An ascending sort enables the selection of the $q$ most diverse acquisitions as the final $q$ elements of this sorted vector.

 \subsection{Nystr\"{o}m approximation for scalability in high dimensions}
\label{ss:nystrom}
In high dimensions, the posterior sampling in \qpots~becomes computationally expensive. We address this via a Nystr\"{o}m approximation~\cite{williamsUsingNystromMethod2000} with GP posteriors.  
We begin by writing the reparametrization technique for a GP posterior:
\[
\begin{split}
Y(X) =& \mu(X) + \left[ k(X, X_n) k(X_n, X_n)^{-1} k(X_n, X) \right]^{1/2} Z \\
=& \mu(X) + \Sigma(X)^{1/2} Z,
\end{split}
\]
where $\Sigma(X) = k(X, X_n) k(X_n, X_n)^{-1} k(X_n, X)$ is an $N \times N$ posterior covariance matrix, that we also write as $\Sigma_{NN}$; the main computational bottleneck is in computing the square-root of $\Sigma(X)$ which costs $\mcl{O}(N^3)$. We now choose a subset $X_m \subset X$, $m \ll N$, to write the following, augmented $(N+m) \times (N+m)$ matrix
\[\Sigma_{(N+m)(N+m)} = \begin{bmatrix}
    \Sigma_{mm} & \Sigma_{mN} \\
    \Sigma_{Nm} & \Sigma_{NN}
\end{bmatrix},\]
where $\Sigma_{mm}$ is the intersection of the $m$ rows and $m$ columns of $\Sigma_{NN}$ corresponding to $X_m$, and $\Sigma_{mN} = \Sigma_{Nm}^\top$ are the $m$ rows of $\Sigma_{NN}$ corresponding to $X_m$. Then, we have that
\begin{equation}
    \begin{split}            
    \Sigma(X) =& \Sigma_{NN} \approx \Sigma_{mN} \Sigma_{mm}^{-1} \Sigma_{mN}^\top\\
    =& \left(\Sigma_{mN} \Sigma_{mm}^{-1/2} \right)\left(\Sigma_{mN} \Sigma_{mm}^{-1/2} \right)^\top.
    \end{split}
    \label{eqn:nystrom_squareroot}
\end{equation}
Therefore, we have the Nystr\"{o}m approximation of the posterior covariance matrix square root as $\Sigma(X)^{1/2} \approx  \left(\Sigma_{mN} \Sigma_{mm}^{-1/2} \right).$
The approximation in \Cref{eqn:nystrom_squareroot} requires an $m\times m$ matrix square root which costs $\mcl{O}(m^3)$ and an additional $\mcl{O}(Nm^2)$ for the matrix multiplication. Therefore, an overall cost of $\mcl{O}(m^3 + Nm^2)$ is much better than $\mcl{O}(N^3)$, provided $m \ll N$. The choice of $X_m$ is an open research question; in this work, we choose them to be the subset of nondominated points in $\{X_n, \mbf{y}_n\}$ which works uniformly well in practice.

In summary, the core idea of our method is inspired by Thompson sampling in single-objective optimization: while TS selects candidate points with high probability of being optimal, our approach selects points based on their probability of being Pareto optimal. \qpots~ generates acquisitions by solving a computationally inexpensive multiobjective optimization problem on posterior GP sample paths. The resulting acquisitions are selected from the predicted Pareto set such that they maximize the minimum distance to previously observed points, thereby promoting diversity and balancing exploration with exploitation. A key advantage of our approach is that selecting a batch of $q$ points incurs no additional cost compared to sequential acquisitions. Simplicity and efficiency are the central strengths of our method.
Our overall methodology is summarized in \Cref{a:method}, and our software implementation is publicly available at \url{https://github.com/csdlpsu/qpots}.

\begin{algorithm}[tb]
 \caption{\qpots: Batch Pareto optimal Thompson sampling}
 \label{a:method}
\begin{algorithmic}
   \STATE {\bfseries Input:} With data sets $\D_n^{1:K}$, fit $K$ GP models $\mcl{M}_1$ through $\mcl{M}_{K}$
     and GP hyperparameters $\bs{\Omega}_1,\ldots,\bs{\Omega}_{K}$. \\
     Parameters $B$ (total budget) \\
   \STATE \textbf{Output:} {Pareto optimal solution $\mcl{X}^*,~\mcl{Y}^*$}.
  \STATE \FOR{$i=n+1, \ldots, B$, }
      \STATE {\bf 1.} Sample the posterior GP sample paths: $Y_i(\x, \omega),~\forall i=1,\ldots,K$.
      \STATE  {\bf 2.} {\bf (Optional)}\ Use Nystr\"{o}m approximation for scalability in higher dimensions.
      \STATE {\bf 3.} Solve the cheap multiobjective optimization problem  \Cref{eqn:moo_gp} via evolutionary algorithms. 
      \STATE {\bf 4.} Choose $q$ candidate points from $X^*$ according to \eqref{eqn:maximin}.
      \STATE {\bf 5.} Evaluate objectives and append to data set $\D^{1:K}_i$.
      \STATE {\bf 6.} Update GP hyperparameters $\bs{\Omega}_1,\ldots,\bs{\Omega}_{K}$
    \ENDFOR
\end{algorithmic}
\end{algorithm}

 \subsection{Related work}
\label{sec:related_work}

Naturally, one may seek to extend ideas from single-objective Bayesian optimization (BO) to the multiobjective setting; indeed, this intuition underlies much of the development in multiobjective Bayesian optimization (MOBO). Among the most prevalent approaches are those based on \emph{scalarization}, wherein the $K$ objectives are combined into a single scalar objective, allowing the application of standard single-objective BO methods. For instance, multiobjective efficient global optimization (ParEGO)~\cite{knowles2006parego} performs scalarization using randomly drawn weights and applies expected improvement (EI) to the resulting objective. Multiobjective evolutionary algorithm based on decomposition (MOEA/D-EGO)~\cite{zhang2009expensive} extends ParEGO to batch settings by employing multiple scalarizations and optimizing them in parallel via a genetic algorithm~\cite{zhou2012multiobjective}. Other extensions include qParEGO and its noisy counterpart, qNParEGO~\cite{daultonParallelBayesianOptimization2021}, designed for batch and noisy evaluations, respectively.
Paria et al.~\cite{paria2020flexible} similarly apply Chebyshev scalarization but use Thompson sampling (TS) to optimize the scalarized objective. Zhang et al.~\cite{zhang2020random} propose a hypervolume-based scalarization whose expected value, under a specific distribution over weights, is equivalent to the hypervolume indicator. Their method adopts an upper confidence bound framework, although it is evaluated only on limited benchmarks. Despite their simplicity, scalarization-based approaches tend to be sample inefficient and often struggle with non-convex, disconnected, or otherwise complex Pareto frontiers.

Beyond scalarization, many MOBO methods directly construct acquisition functions based on the \emph{hypervolume indicator}. These approaches aim to quantify expected improvement in hypervolume, but doing so typically leads to non-convex optimization problems that are difficult to solve. The expected hypervolume improvement (EHVI)~\cite{couckuytFastCalculationMultiobjective2014, emmerichComputationExpectedImprovement2008} is a canonical extension of EI to the multiobjective setting, with several gradient-based variants also available~\cite{yang2019multi, daulton2020differentiable}. Likewise, the stepwise uncertainty reduction (SUR) criterion has been adapted for MOBO~\cite{picheny2015multiobjective}. More recent methods~\cite{daultonDifferentiableExpectedHypervolume2020, daultonParallelBayesianOptimization2021} extend EHVI to batch acquisitions ($q > 1$) with differentiable formulations to facilitate gradient-based optimization. Nevertheless, these approaches often suffer from noisy acquisition functions and increasing optimization difficulty as batch size grows.

Recognizing the challenges posed by the non-convexity of hypervolume-based acquisition functions, several recent methods introduce heuristics to better balance exploration and exploitation. Diversity-Guided Efficient Multiobjective Optimization (DGEMO)~\cite{konakoviclukovicDiversityGuidedMultiObjectiveBayesian2020} selects diverse candidates by solving local optimization problems that improve hypervolume, and scales well to large batch sizes. However, DGEMO does not handle noisy observations. Pareto Active Learning (PAL)~\cite{zuluaga2013active} takes a classification-based approach by assigning to each candidate the probability of being Pareto optimal, based on GP posterior uncertainty. The resulting acquisition function is again difficult to optimize and is typically addressed by discretizing the design space~\cite{zuluaga2013active}.

An alternative to hypervolume improvement is offered by \emph{entropy-based methods}, which quantifies expected reduction in uncertainty of the Pareto frontier. Predictive entropy search for multiobjective optimization (PESMO)~\cite{hernandez-lobatoPredictiveEntropySearch} and Max-value Entropy Search for Multiobjective optimization (MESMO)~\cite{belakariaMaxvalueEntropySearch2019} are two prominent examples. MESMO improves on PESMO in scalability and computational cost, but still lacks a closed-form acquisition function and depends on sampling approximations. Pareto Frontier Entropy Search (PFES)~\cite{suzuki2020multi} offers a closed-form entropy criterion via random Fourier features (RFFs), but is restricted to stationary GP kernels. Similarly, Joint Entropy Search for Multiobjective Optimization (JESMO)~\cite{hvarfner2022joint, tu2022joint} estimates joint information gain over inputs and outputs, but still requires approximations to the acquisition function.

Despite these developments, the use of TS in MOBO remains relatively unexplored. Existing methods that incorporate TS often do so in conjunction with traditional scalarization or hypervolume techniques, inheriting their associated limitations. For example, Thompson Sampling Efficient Multiobjective Optimization (TSEMO)~\cite{bradford2018efficient} samples Pareto sets from GP posteriors and evaluates them using hypervolume improvement. TS with Chebyshev scalarization (TS-TCH) likewise depends on scalarized objectives. Multiobjective Regionalized Bayesian Optimization (MORBO)~\cite{daulton2022multi} introduces a trust-region strategy and uses TS to identify local candidates; however, it ultimately relies on hypervolume improvement for selection.

\paragraph{Existing limitations and our contributions.}  
In summary, current MOBO methods exhibit one or more of the following limitations:  
(1) Reliance on acquisition functions that lead to difficult inner optimization problems, reducing sample efficiency;   
(2) Limited or no support for batch acquisitions; and  
(3) Difficulty handling noisy objectives. Our proposed approach addresses these challenges by targeting three key desiderata:  
(i) sample efficiency in both sequential and batch settings,    
(ii) robustness to noisy objective evaluations, and  
(iii) simplicity of implementation. 

\mg{
The advantages of treating multipoint design as a multi-objective problem should be more explicitly articulated: why and when Pareto sets are preferable to weighted-sum scalarization in aerospace contexts.
- Provide a stronger aerodynamic interpretation of the Pareto front in the CRM case to demonstrate the practical value of multiple solutions, not just numerical improvement.
}

\section{Experiments}
\label{sec:experiments}
We now demonstrate our proposed approach, \qpots, on several synthetic and real-world experiments of varying input ($d$) and \mg{keeping output ($K$) dimensions to $K=2$, to keep it consistent with the CRM experiment although our method seamlessly extends to any $K$}. We compare \qpots~against TSEMO (another popular TS-based approach), PESMO, JESMO, MESMO (entropy-based approaches), qNEHVI, qPAREGO (classical approaches), and random (Sobol) sampling.

\begin{table}[]
    \centering
    \begin{tabular}{c|c|c|c|c|c}
         \hline \hline
         {\bf Experiment} & Branin-Currin & DTLZ3 & DTLZ7 & ZDT3 & NASA CRM \\ 
         \hline
         {\bf Reference point} & $(70, 70)$ & $(10^4, 10^4)$ & $(25, 25)$ & $(11, 11)$ & $(10^2, 10^2)$ \\
        \hline
        {\bf Dimensionality $(d, K)$} & $(2,2)$ & $(10,2)$ & $(5,2)$ & $(10,2)$ & $(24, 2)$ \\
        \hline
        {\bf Seed/repetitions} & $20/10$ & $100/10$ & $50/10$ & $100/10$ & $200/1$ \\
        \hline
        {\bf Batch size ($q$)} & $1,4$ & $1,4$ & $1,4$ & $1,4$ & $1$ \\
        \hline
    \end{tabular}
    \caption{Summary of experiment settings.}
    \label{tab:ref_point}
\end{table}
\subsection{Synthetic experiments}
\label{sec:synthetic}
We consider the following synthetic experiments: the Branin-Currin function, the ZDT3 function, which has disconnected Pareto frontiers, the DTLZ3 function, and the DTLZ7 test function. 
We provide each experiment with a set of $n=10\times d$ seed samples and observations to start the algorithm, chosen uniformly at random from $\mcl{X}$; then they are repeated $10$ times. However, the seed is fixed for each repetition across all acquisition policies to keep the comparison fair.  We add Gaussian noise with variance $\tau^2=10^{-3}$ to all our observations to assist with the conditioning of the GP covariance matrix.
 The experiments are run with a high-performance compute node with 256GB memory and $48$-core nodes; repetitions are parallelized across the cores. 
Our metric for comparison is the hypervolume indicator computed via box decomposition~\cite{yang2019efficient}. \mg{It should be noted that hypervolume computation is a postprocessing step in the methodology and hence is not part of the algorithm. The hypervolume computation involves specifying a reference point (which is the absolute worst value each objective can take) -- these are summarized in \Cref{tab:ref_point}. For synthetic experiment, we determine this by densely sampling the input domain $\X$ and evaluating the objectives on these samples. For the NASA CRM experiment, we use domain knowledge (that drag coefficients are lower bounded by $0$ and hence set $(100,100)$. }

Throughout the manuscript, we default to the anisotropic Mat\'{e}rn class kernel with $\nu=5/2$.  Our implementation primarily builds on GPyTorch~\cite{gardner2018gpytorch} and BoTorch~\cite{balandat2019botorch}; we leverage 
MPI4Py~\cite{dalcin2021mpi4py} to parallelize our repetitions. We use PyMOO~\cite{blank2020pymoo} to leverage its NSGA-II solver in our inner optimization (\Cref{eqn:moo_gp}); the population size for all problems is fixed at $100\times d$. 

\begin{figure}[htb!]
    \centering
    \begin{subfigure}{.5\textwidth}
        \includegraphics[width=.95\linewidth]{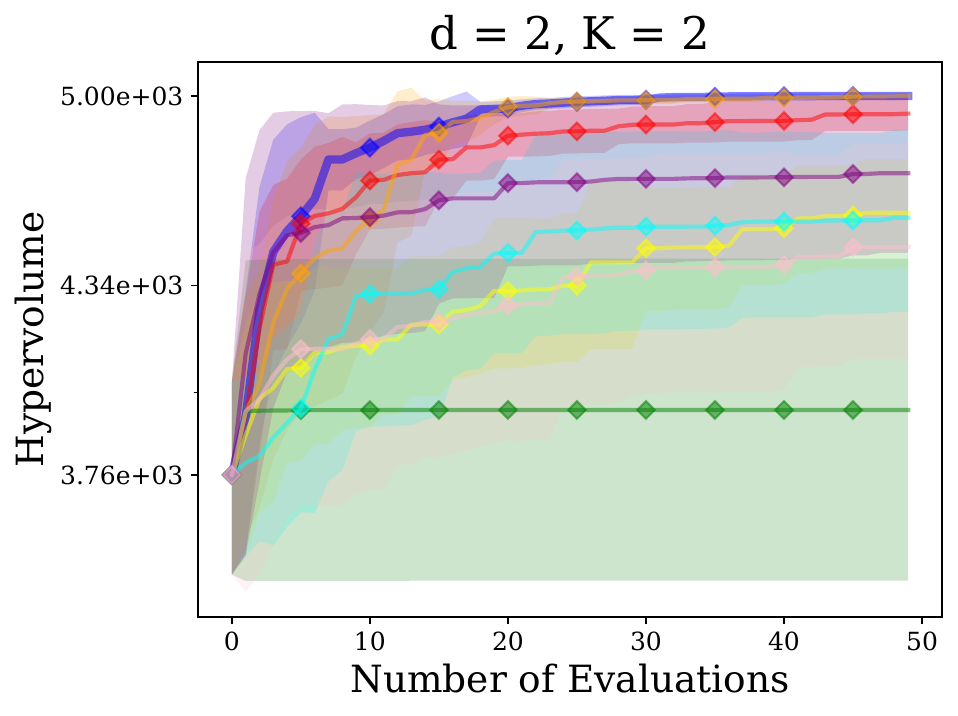}
        \caption{Branin-Currin}
    \end{subfigure}%
    \begin{subfigure}{.5\textwidth}
        \includegraphics[width=.95\linewidth]{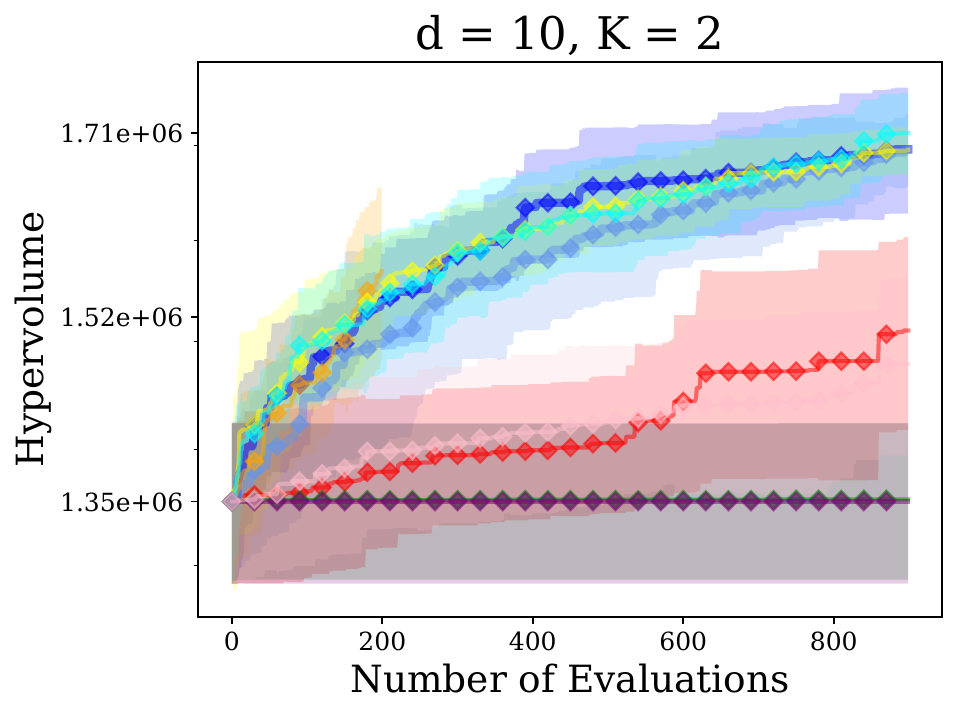}
        \caption{DTLZ3}
    \end{subfigure}\\
    \begin{subfigure}{.5\textwidth}
        \includegraphics[width=.95\linewidth]{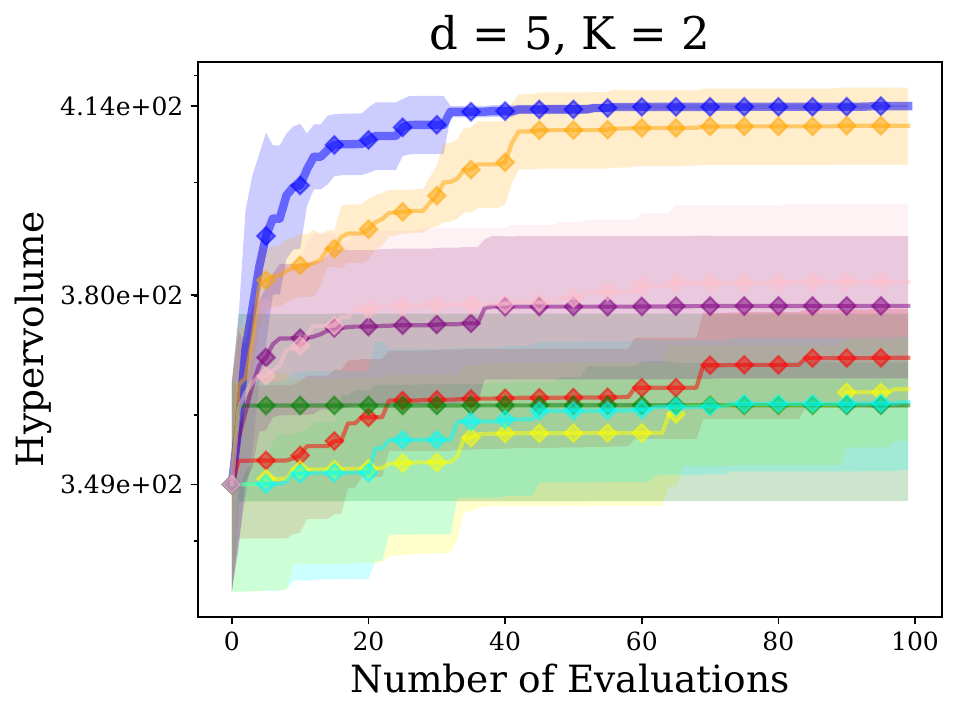}
        \caption{DTLZ7}
    \end{subfigure}%
    \begin{subfigure}{.5\textwidth}
        \includegraphics[width=.95\linewidth]{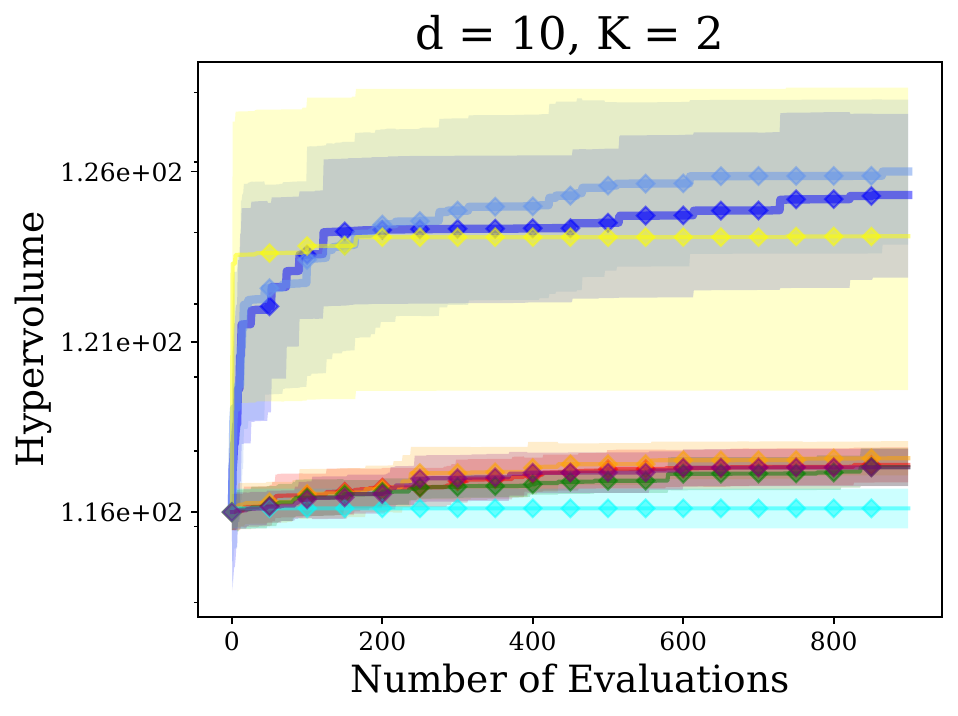}
        \caption{ZDT3}
    \end{subfigure}\\
    \begin{subfigure}{1\textwidth}
        \includegraphics[width=1\linewidth]{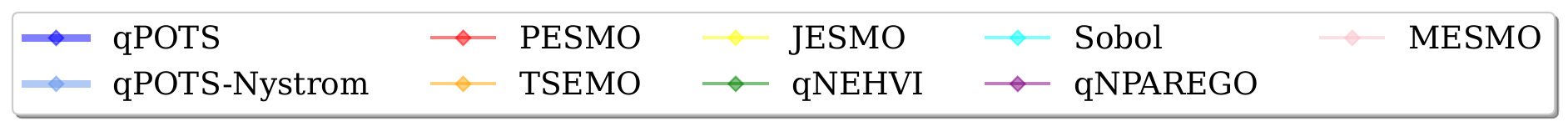}
    \end{subfigure}\\
    \caption{Hypervolume (higher the better) Vs. iterations for sequential ($q=1$) acquisition; plots show mean and $\pm 1$ standard deviation out of $10$ repetitions.}
    \label{fig:hv_seq}
\end{figure}

\Cref{fig:hv_seq} shows the performance comparison (in terms of hypervolume) for sequential sampling ($q=1$). Notice that \qpots~(blue) is consistently the best
performer, or is amongst the best. We show
\qpots-Nystrom-Pareto (that is $X_m$ chosen as the current Pareto set) for the $d=10$ experiments which shows no loss of accuracy for the DTLZ3 function and, surprisingly, a slight gain in accuracy for the ZDT3 function. Despite the success in the $q=1$ setting, our main benefit is seen in the batch sampling setting which we discuss next.

\begin{figure}[htb!]
    \centering
    \begin{subfigure}{.5\textwidth}
        \includegraphics[width=.95\linewidth]{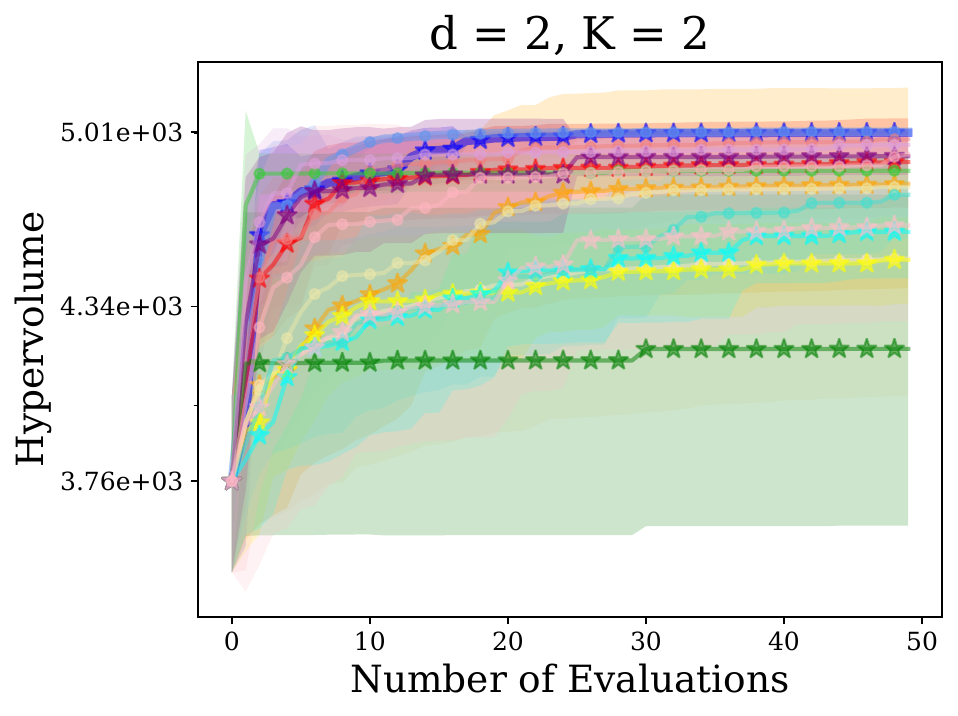}
        \caption{Branin-Currin}
    \end{subfigure}%
    \begin{subfigure}{.5\textwidth}
        \includegraphics[width=.95\linewidth]{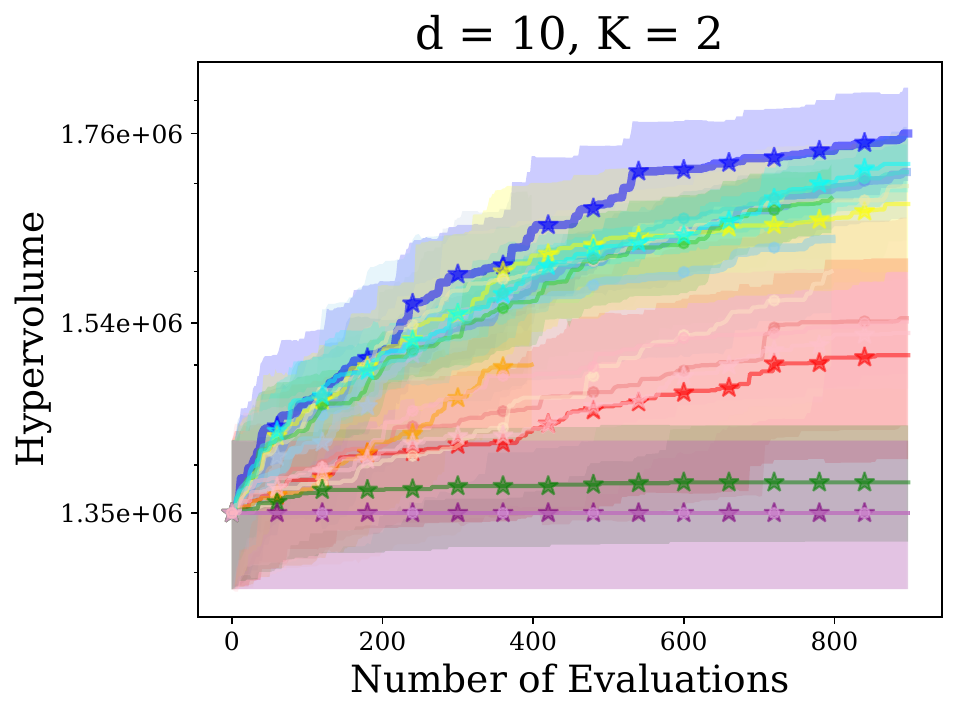}
        \caption{DTLZ3}
    \end{subfigure}\\
    \begin{subfigure}{.5\textwidth}
        \includegraphics[width=.95\linewidth]{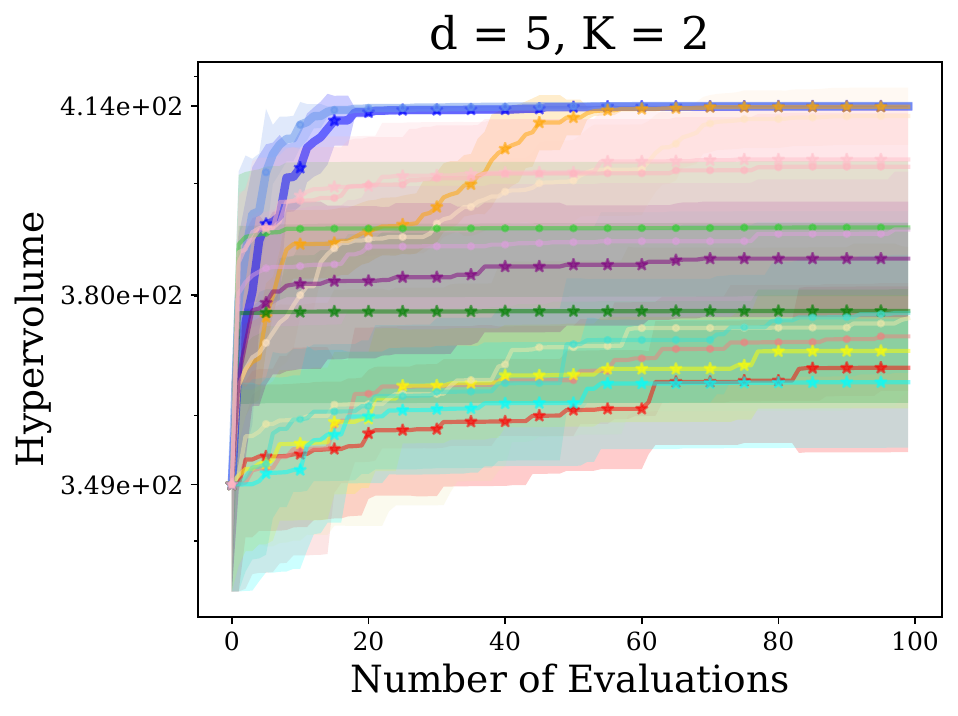}
        \caption{DTLZ7}
    \end{subfigure}%
    \begin{subfigure}{.5\textwidth}
        \includegraphics[width=.95\linewidth]{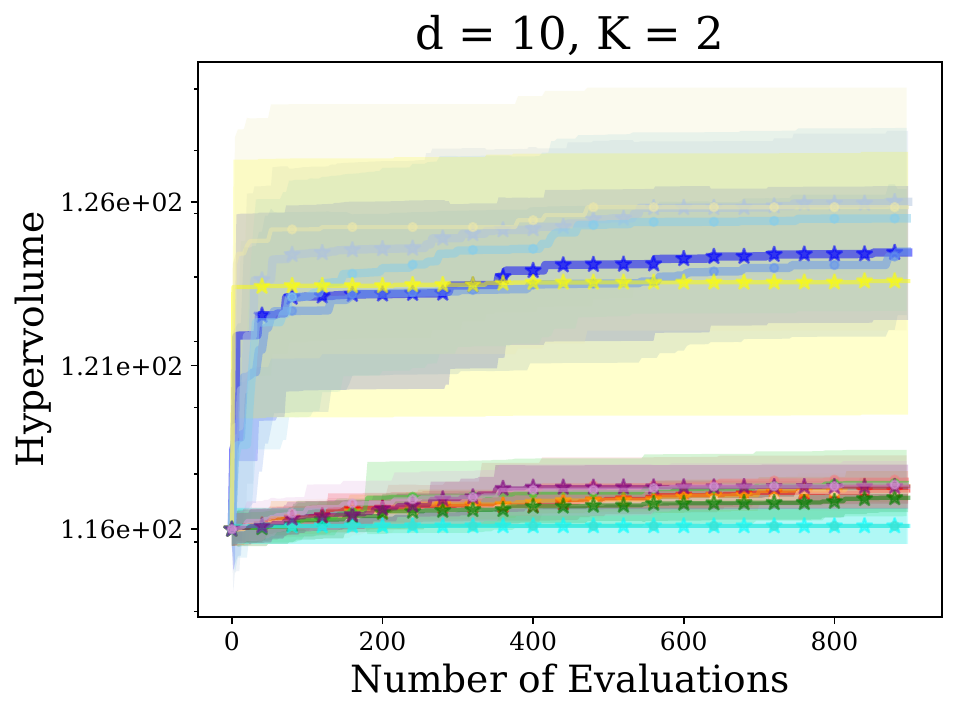}
        \caption{ZDT3}
    \end{subfigure}\\
    \begin{subfigure}{1\textwidth}
        \includegraphics[width=1\linewidth]{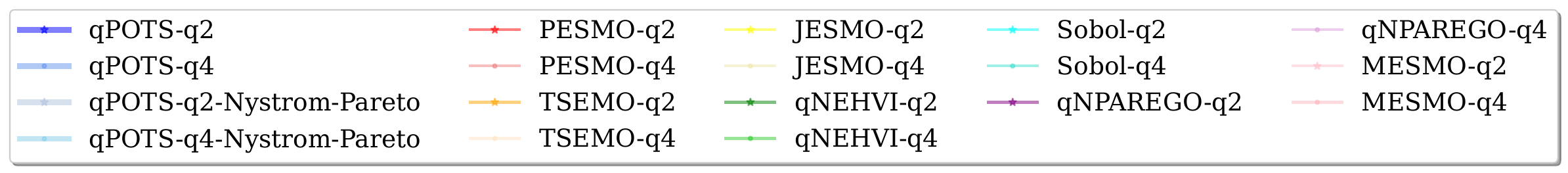}
    \end{subfigure}\\
    \caption{{\bf Batch acquisition.} Hypervolume Vs. iterations for batch ($q>1$) acquisition; plots show mean and $\pm 1$ standard deviation out of $10$ repetitions. \qpots~outperforms all competitors, but the benefit is more pronounced in the batch case. }
    \label{fig:hv_batch}
\end{figure}

The demonstration for the batch setting is shown in \Cref{fig:hv_batch}, where we present experiments with $q=2, 4$; like in the $q=1$ setting we include \qpots-Nystrom-Pareto for the $d=10$ experiments only.  While \qpots~comfortably wins in all experiments, notice that the difference with TSEMO is more pronounced in the batch setting. The advantage in the batch setting $q>1$ is a key benefit of \qpots, because the cost of batch sampling is practically the same as sequential ($q=1$) sampling. \mg{It is expected that \qpots~performs better in the batch setting; this is for the following reasons. While the Thompson sampling automatically balances exploration and exploitation in the iterates, we further ensure diversity in the candidates by sorting them according to the maximin distance. Additionally, our batch sampling avoids complex stochastic optimization as in other methods (e.g., qNEHVI) which makes the computation cheap and deterministic.} 

\subsection{The NASA CRM experiment}
\label{sec:crm}
\begin{figure}[htb!]
    \centering
    \begin{subfigure}{0.5\textwidth} 
    \centering
    \includegraphics[width=.9\linewidth,height=6cm]{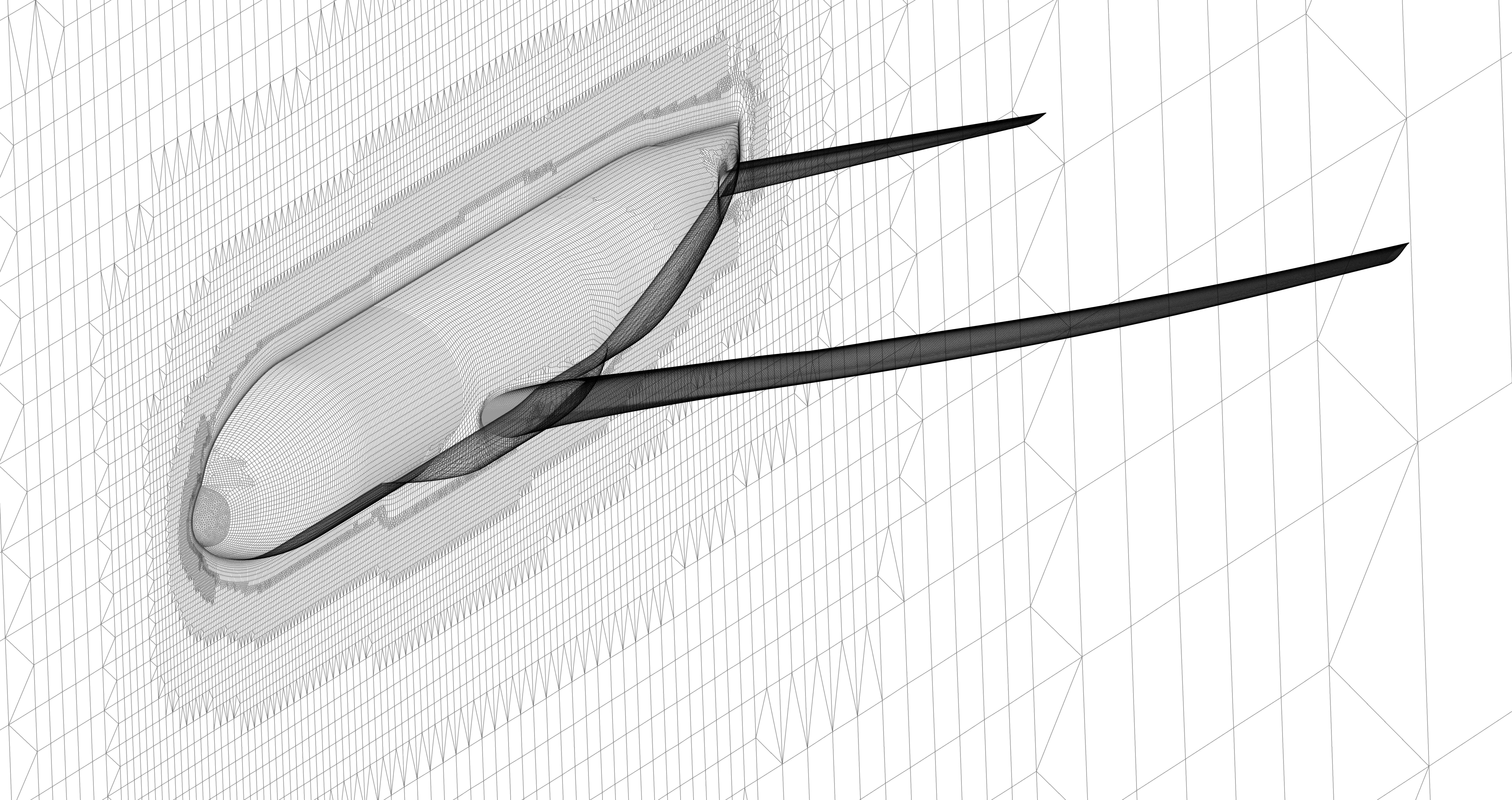}     
    \end{subfigure}%
    \begin{subfigure}{0.5\textwidth}
        \centering
        \includegraphics[width=.9\linewidth, height=6cm]{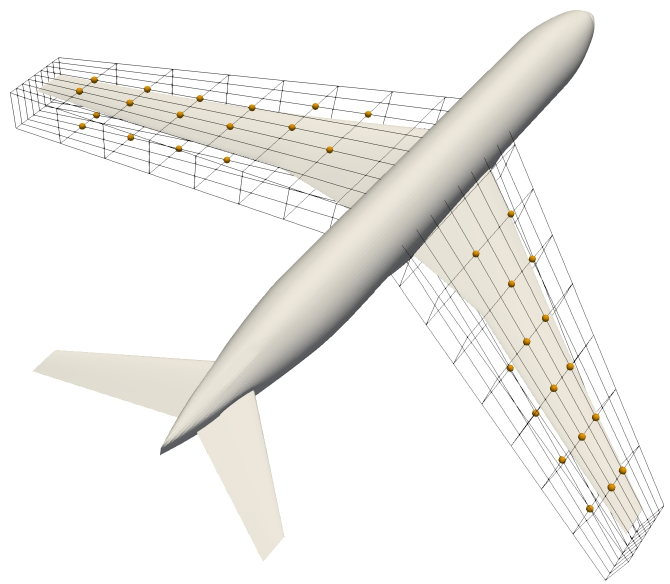} 
    \end{subfigure}
    \caption{Left: Computational mesh of the NASA CRM. Right:The FFD box parameterization of the NASA CRM. The orange dots represent the control points used for the deformation. Only one half of the aircraft was used for the CFD because of symmetry.}
    \label{fig:mesh_and_ffd}
\end{figure}

The NASA common research model (CRM) \citep{vassberg2008development} is an aircraft geometry made for industry-relevant CFD research. We consider a ``clean wing'' (that is, without nacelles and pylons) version of the CRM geometry, as shown in \Cref{fig:mesh_and_ffd}. We leverage the SU2 \citep{economon2016su2} open-source multiphysics code to perform an inviscid Euler simulation of the flow past the CRM. The computational mesh contains $6.12$M volume elements. A scalar upwind scheme was used for the numerical method with the Venkatakrishnan limiter \citep{venkatakrishnan1995convergence}. A constant CFL number of $10$ was used after it was deemed adequate for satisfactory convergence. Our full SU2 input files and mesh are publicly available at \url{https://github.com/csdlpsu/su2scrips/tree/main/parallel_test}.

We consider a two-objective $K=2$ optimization problem, identifying two cruise points as summarized in \cref{tab:flow-cons}; $C_{D1}$ and $C_{D2}$ denote the corresponding drag coefficients at the two cruise points. The optimization problem we will be solving is
\begin{equation}
    \begin{split}
        \min_{\x \in \mathcal{X}} & \quad \{C_{D1}(\textbf{x}), C_{D2}(\textbf{x})\}, \\
    \end{split}
\end{equation}
which is an unconstrained two-objective optimization ($K=2$) problem in $d=24$ dimensions.

\begin{table}[tb!]
    \centering
    \begin{tabular}{|c|c|c|}
        \hline
         & Cruise point 1 & Cruise point 2 \\
         \hline
        Mach number ($M$) & 0.85 & 0.7 \\
        \hline
        Angle of attack ($\alpha$) & 2.37 & 3.87 \\
        \hline
    \end{tabular}
    \caption{Operating conditions for the two chosen cruise points.}
    \label{tab:flow-cons}
\end{table}

We parameterize the wing with a free-form deformation (FFD)~\citep{sederberg1986free} box with 162 control points, $24$ of which are allowed to deform -- these are the design variables in our optimization problem; \Cref{fig:mesh_and_ffd} displays the FFD box, highlighting the design variables. For control points nearer to the center of the wing, the deformation bounds were (0, 1e-2), slightly farther from the center were (0, 1e-3), and control points near the root and the tip had bounds of (0, 1e-5) to prevent CFD solver divergence. 

The CFD simulations were run on Penn State's Roar Collab supercomputer. We use 6 compute nodes with $48$ cores each for a total of $288$ cores and 2400 GB of memory. For GP training, we first provide $200$ seed points. We run the optimization for an additional $200$ iterations with $q=1$. We compare \qpots~ to qNEHVI, qPAREGO, which were the best-performing approaches from synthetic experiments, and Sobol random sampling using the hypervolume metric. Due to the computational cost of the CRM experiment the experiments are not repeated, and our comparison is limited to only qNEHVI and qPAREGO. \mg{Note that BO is particularly suited for optimization of expensive functions where the goal is sample efficiency; therefore, replicating the real-world experiments is not practical. On the other hand, the synthetic experiments are repeated multiple times to establish statistical significance of the reported results. This is standard practice in the Bayesian optimization community~\cite{Liang2021, Shields2021, morita2022applying}. Finally, we don't turn on the Nystr\"{o}m approximation for the CRM experiment. This is because the acquisition computations were still substatially faster than qNEHVI and qPAREGO in this high-dimensional case. }

\begin{figure}[htb!]
    \centering
    \begin{subfigure}{0.5\textwidth}
        \includegraphics[width=1\linewidth, height=6cm]{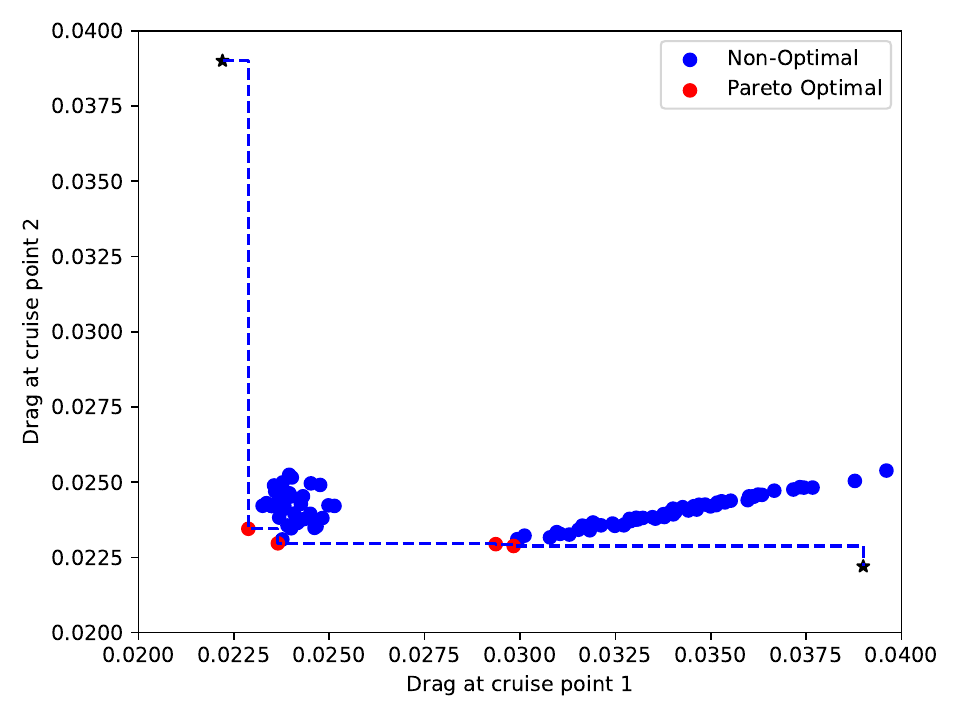}    
    \end{subfigure}%
    \begin{subfigure}{0.5\textwidth}
        \includegraphics[width=1\linewidth, height=6cm]{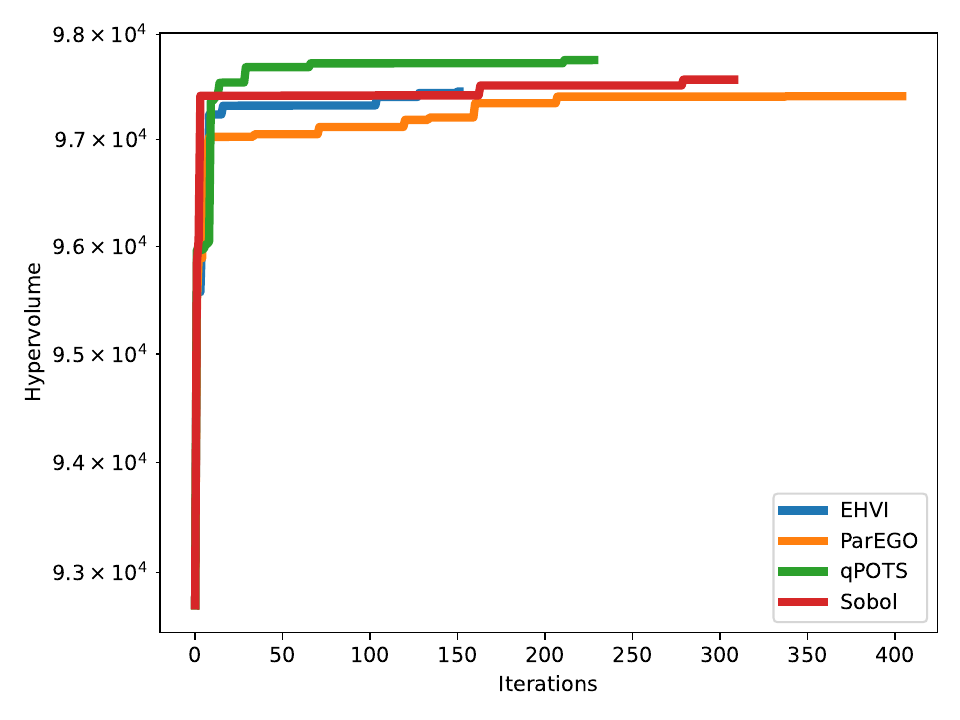}    
    \end{subfigure}
    \caption{Left: Scatter plot of the drag coefficients (blue) and the Pareto optimal values (red) returned by \qpots. Right: Comparison of hypervolume (higher the better) history of \qpots~against that of other methods.}
    \label{fig:hv-comp-aso}
\end{figure}

We begin by presenting the hypervolume history and the Pareto frontiers identified by \qpots, in \Cref{fig:hv-comp-aso}. The left panel shows a scatterplot of the objectives chosen by \qpots~at each iteration. This reveals an important hypothesis made at the beginning of this manuscript -- that the objectives are nonconflicting in parts of the design space while being conflicting in the other parts. Crucially, the conflicting part is where the Pareto optimal designs are also found, emphasizing the importance of multiobjective optimization. The right panel of the figure shows the hypervolume histories of \qpots~alongside competing methods -- notably, \qpots~achieves the highest hypervolume and is the fastest to converge.

\begin{figure}[htb!]
    \centering
    \begin{tabular}{cc}
        \includegraphics[width=0.48\linewidth]{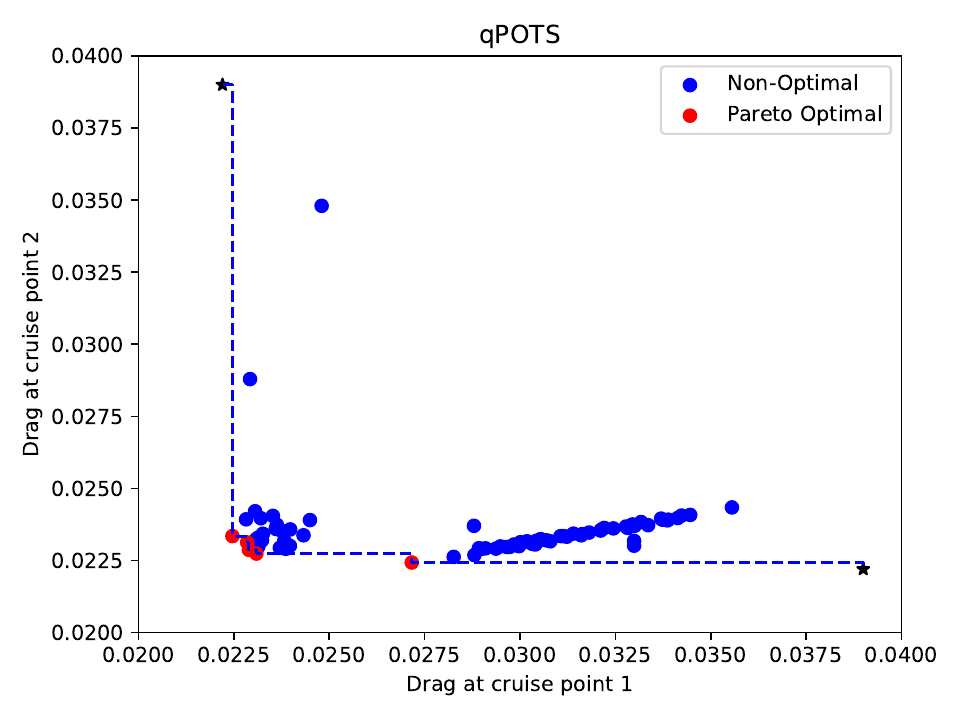} & 
        \includegraphics[width=0.48\linewidth]{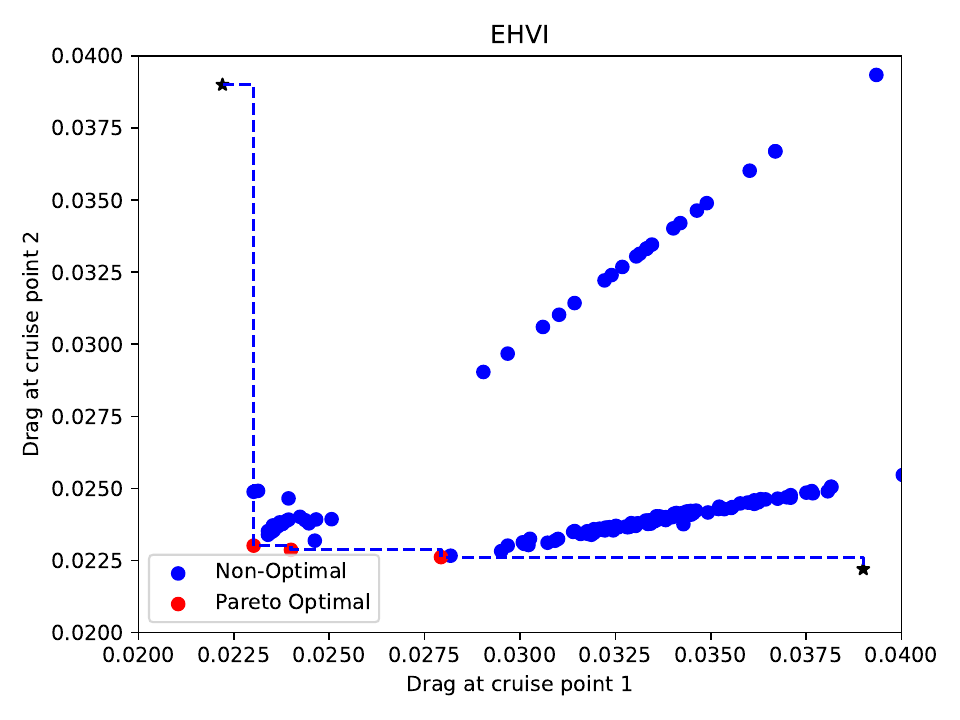} \\
        \includegraphics[width=0.48\linewidth]{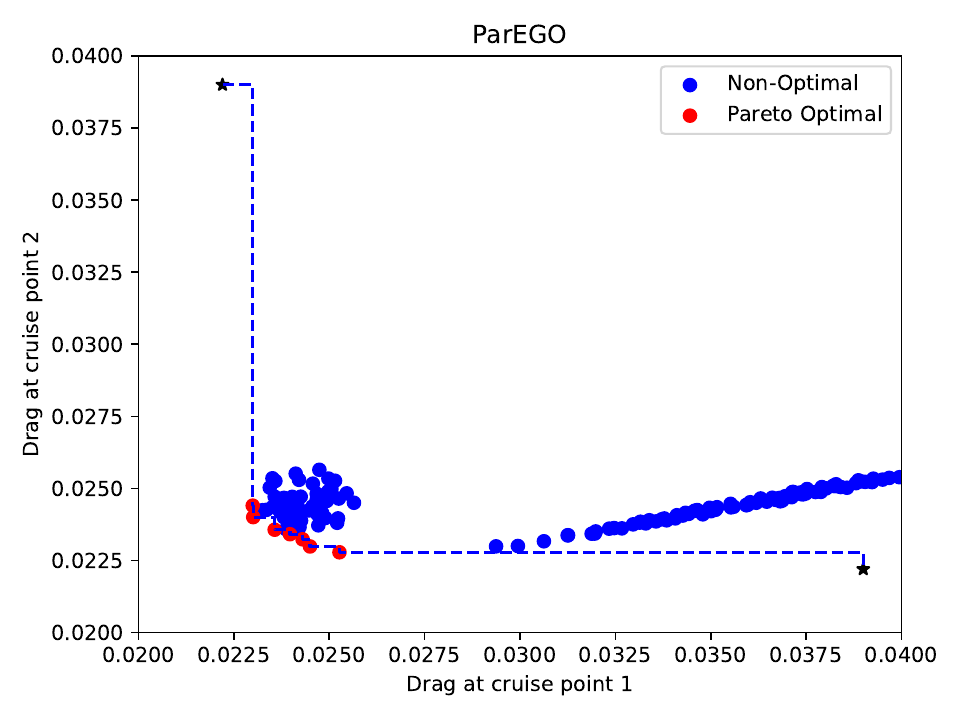} & 
        \includegraphics[width=0.48\linewidth]{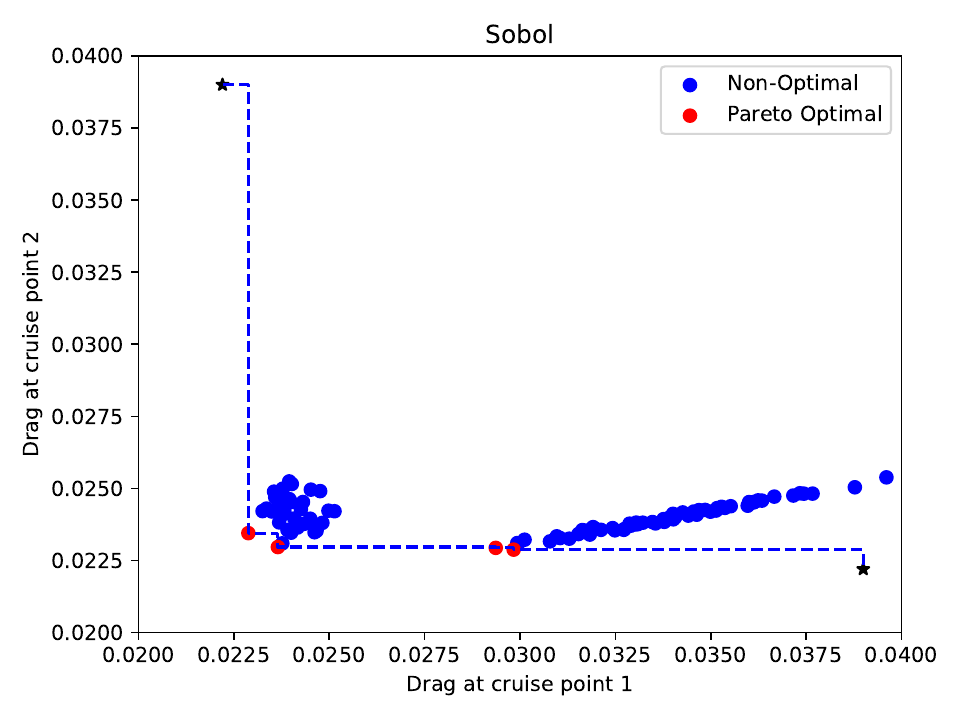} \\
    \end{tabular}
    \caption{Pareto frontiers determined by the various acquisition functions. The blue dots are non-optimal points and the red dots are optimal points. Most of the Pareto optimal points found by the acquisition functions are in the bottom left corner of the graphs, however, there are some Pareto optimal points found with higher drag in the first cruise point.}
    \label{fig:pareto-fronts-crm}
\end{figure}

The scatterplots of objective values are further shown for the other methods compared against, in \Cref{fig:pareto-fronts-crm}. Notice that, similar to \qpots, all methods identify distinct regions of nonconflicting and conflicting objectives. Further, EHVI and Sobol identify only two objective values that are Pareto optimal. While ParEGO, similar to \qpots, identifies several Pareto optimal objectives, the Pareto optimal objectives from ParEGO are dominated by those identified by \qpots. In other words, \qpots~is able to find designs with higher hypervolume compared to other ParEGO (and the rest of the methods) -- an observation that is also confirmed from the right panel of \Cref{fig:hv-comp-aso}. \mg{The aerodynamic interpretation of the Pareto optimal designs are better obtained from the pressure distribution and wing shape plots discussed next.}
\begin{figure}[h!]
    \centering
    \begin{subfigure}{.9\textwidth}
    \centering
        \includegraphics[width=.8\linewidth]{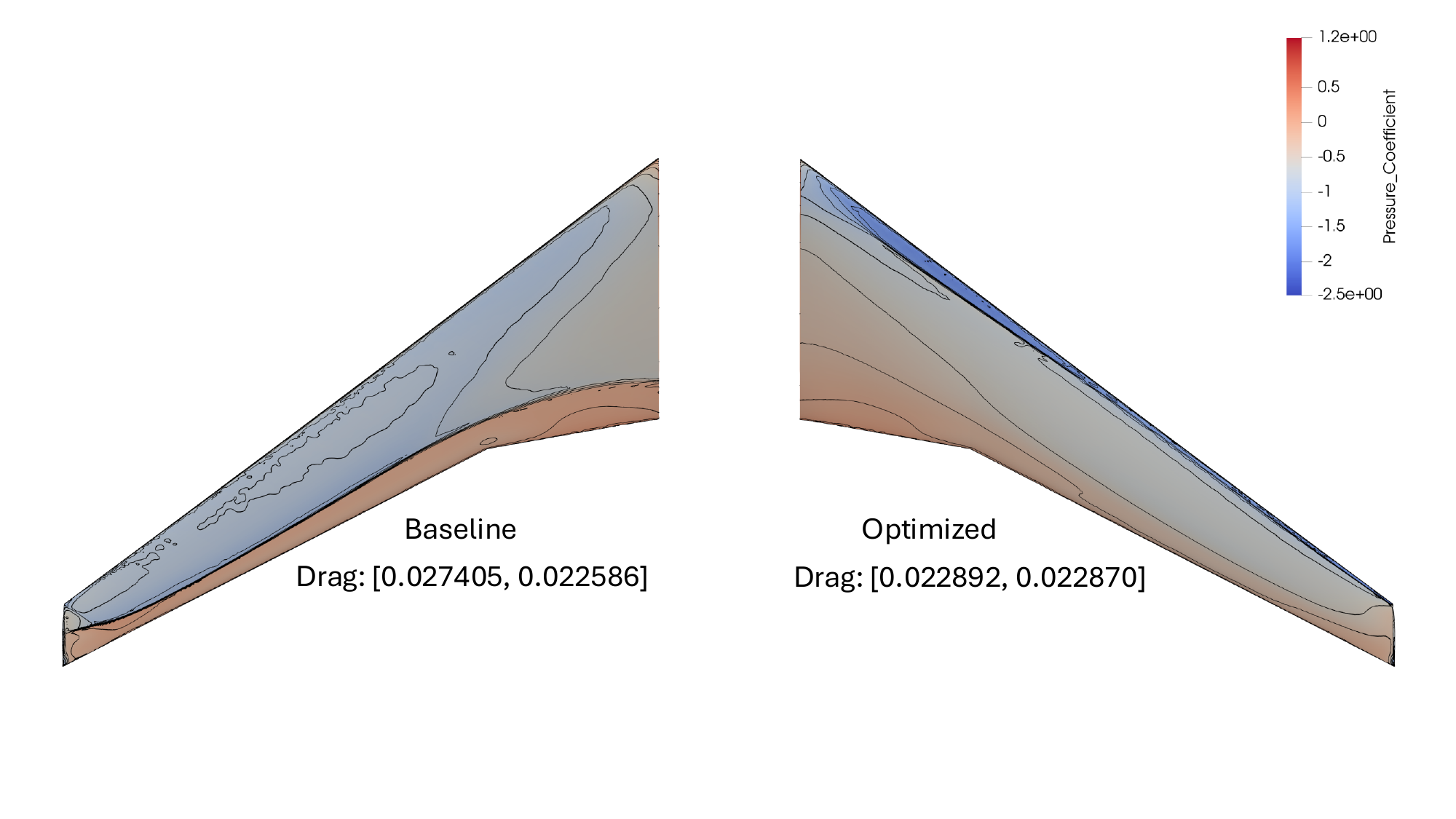}
        \caption{Pareto optimal design $1$.}
    \end{subfigure}\\
    \begin{subfigure}{.9\textwidth}
    \centering
        \includegraphics[width=.8\linewidth]{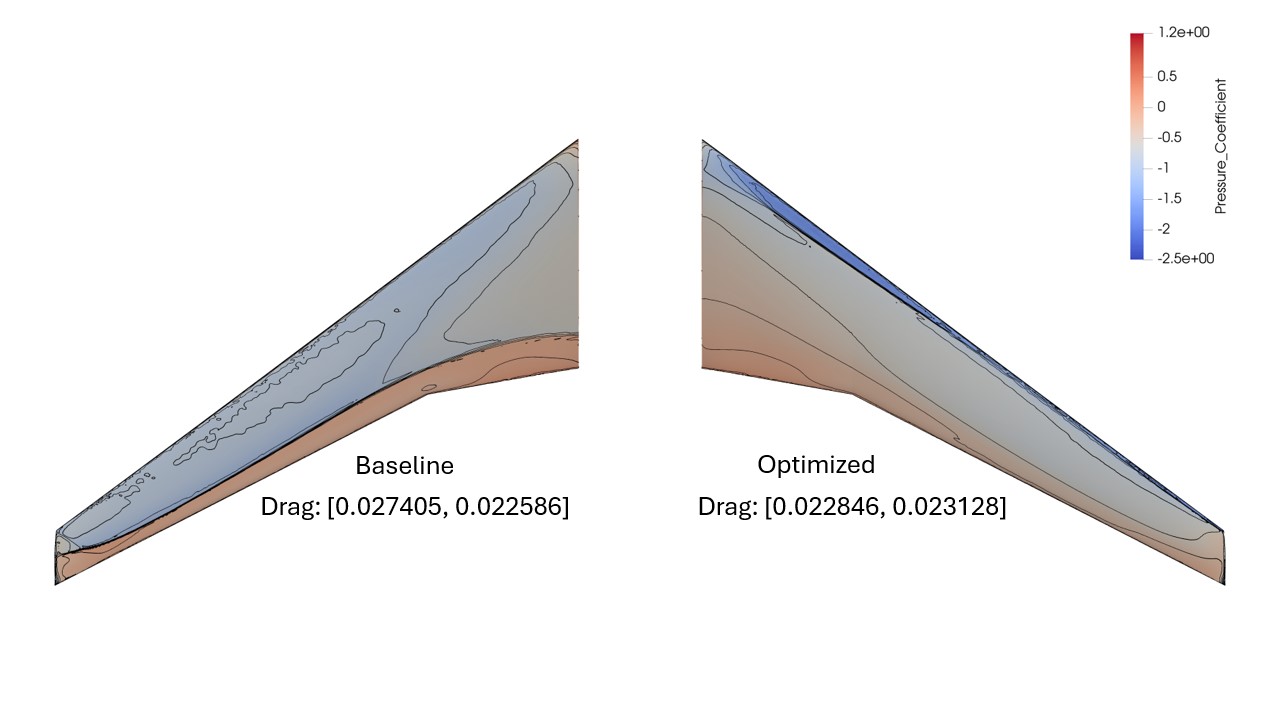}
        \caption{Pareto optimal design $2$.}
    \end{subfigure}\\
    \begin{subfigure}{.9\textwidth}
    \centering
        \includegraphics[width=.8\linewidth]{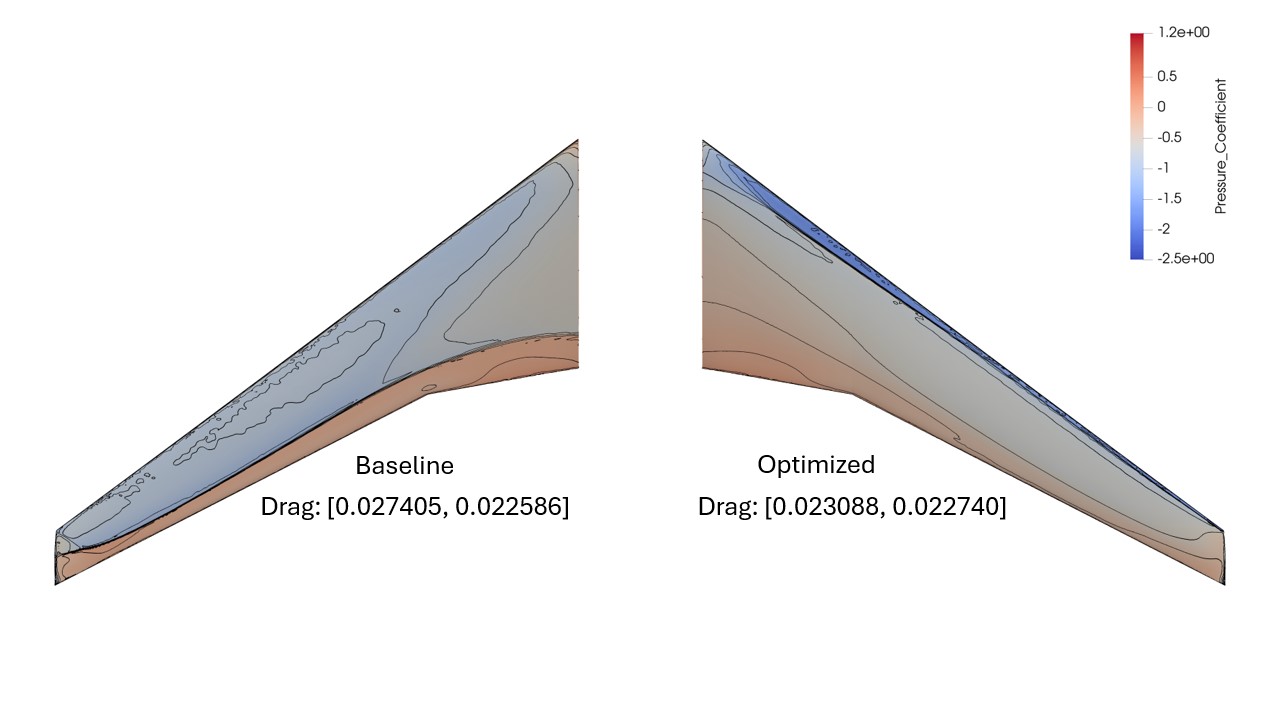}
        \caption{Pareto optimal design $3$.}
    \end{subfigure}
    \caption{Pressure coefficient comparisons for the undeformed baseline wing and two other Pareto optimal points. Again, there is considerable reduced drag in the first objective due to the location change of the shock.}
    \label{fig:wing-comps}
\end{figure}

Next, we compare the pressure coefficient distribution on the wings and  the wing cross-sectional shapes identified by \qpots. We compare them against an undeformed baseline wing to emphasize the impact of optimization; note that this is the unmodified CRM wing without any optimization involved.
\Cref{fig:wing-comps} shows the pressure coefficients of the undeformed baseline wing on the left and an optimized wing at a Pareto optimal point on the right, for three different Pareto optimal designs; note that these correspond to cruise point 1. The drag coefficients at the two cruise points are printed below the wing as $[C_{D,1}, C_{D,2}]$. \mg{First, notice that there is a strong shock toward the trailing edge of the wing for the baseline design. A common thread across all three Pareto optimal designs is that the optimized wing, for cruise point $1$, has a \emph{lower drag by $45-50$ counts}. On the other hand, the optimized drag at cruise point 2 is practically the same (or marginally higher) than the baseline drag for cruise point 2. This shows that the designs at the two cruise points are conflicting -- in other words, an independent optimization at one of the design points is likely to result in an undesirably large drag for the other cruise point, thereby emphasizing the need for multiobjective optimization.}

\begin{figure}[h!]
    \centering
    \includegraphics[width=1\linewidth]{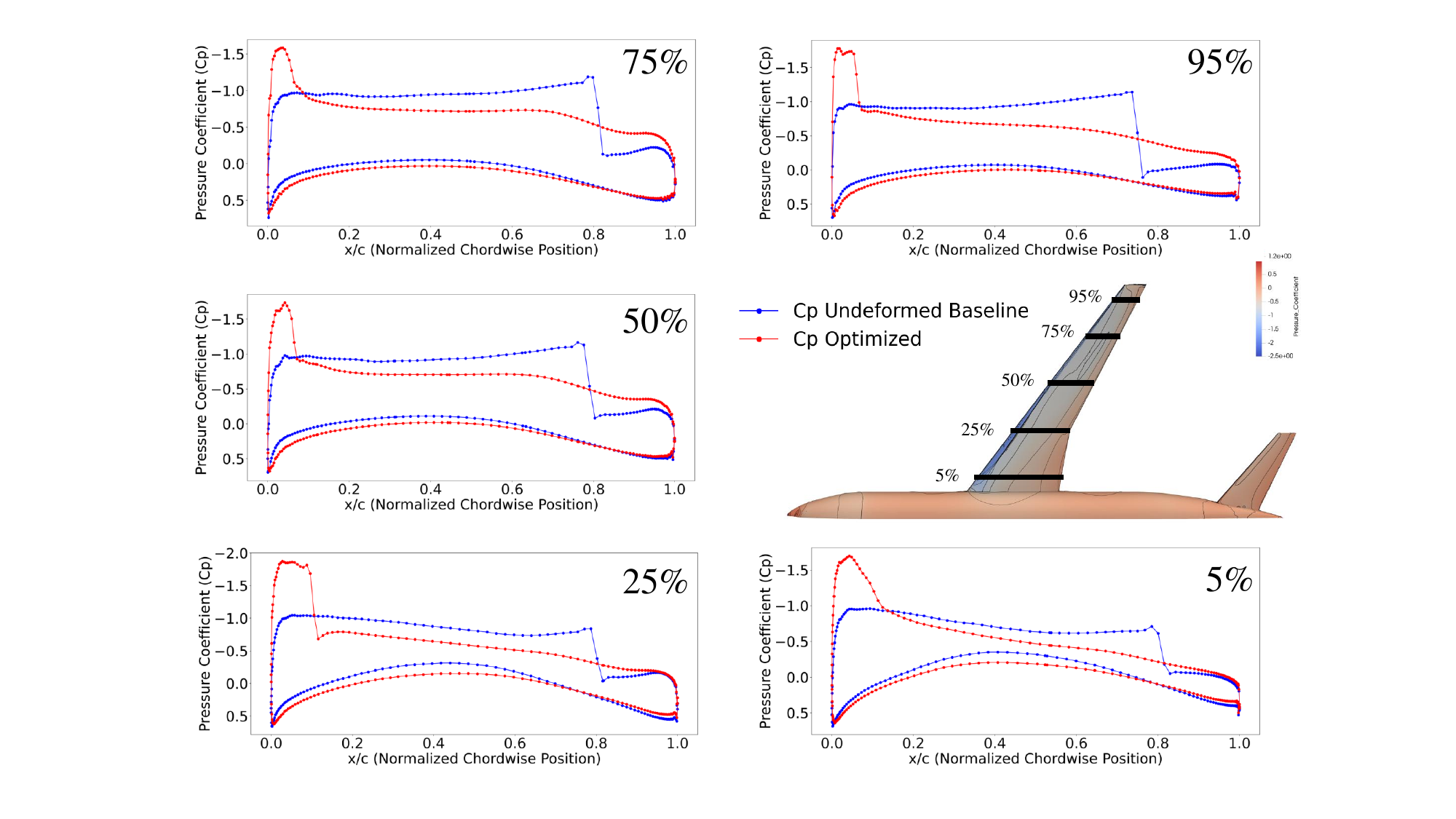}
    \caption{Comparison of pressure coefficients at different percentages of the span for one Pareto optimal point. In blue is the undeformed baseline wing and in red is the optimized wing. Notice the large pressure drops from the undeformed wing around 80\% of the chord, contributing to high drag.}
    \label{fig:cp-wing-sections}
\end{figure}

We next compare the pressure coefficient at selected spanwise cross-sections in \Cref{fig:cp-wing-sections}; note that this figure shows only the Pareto optimal shapes predicted by \qpots~since we already established that \qpots~outperforms other methods. In this plot, the red and blue lines correspond to the optimized and baseline wings, respectively. The relocation of the shock (from the trailing edge to leading edge of the wing) and its reduction in strength is quite evident from these plots. Furthermore, one can appreciate the overall reduction in the wave drag component by visualizing the area under these pressure coefficient plots. Crucially, unlike what could have happened in single objective optimization, multiobjective optimization does not entirely eliminate the shock. This is due to the tradeoffs that must be respected at the other operating condition which, thus, again emphasizes the need for multiobjective optimization for aerodynamic design.

\begin{figure}[h!]
    \centering
    \includegraphics[width=1\linewidth]{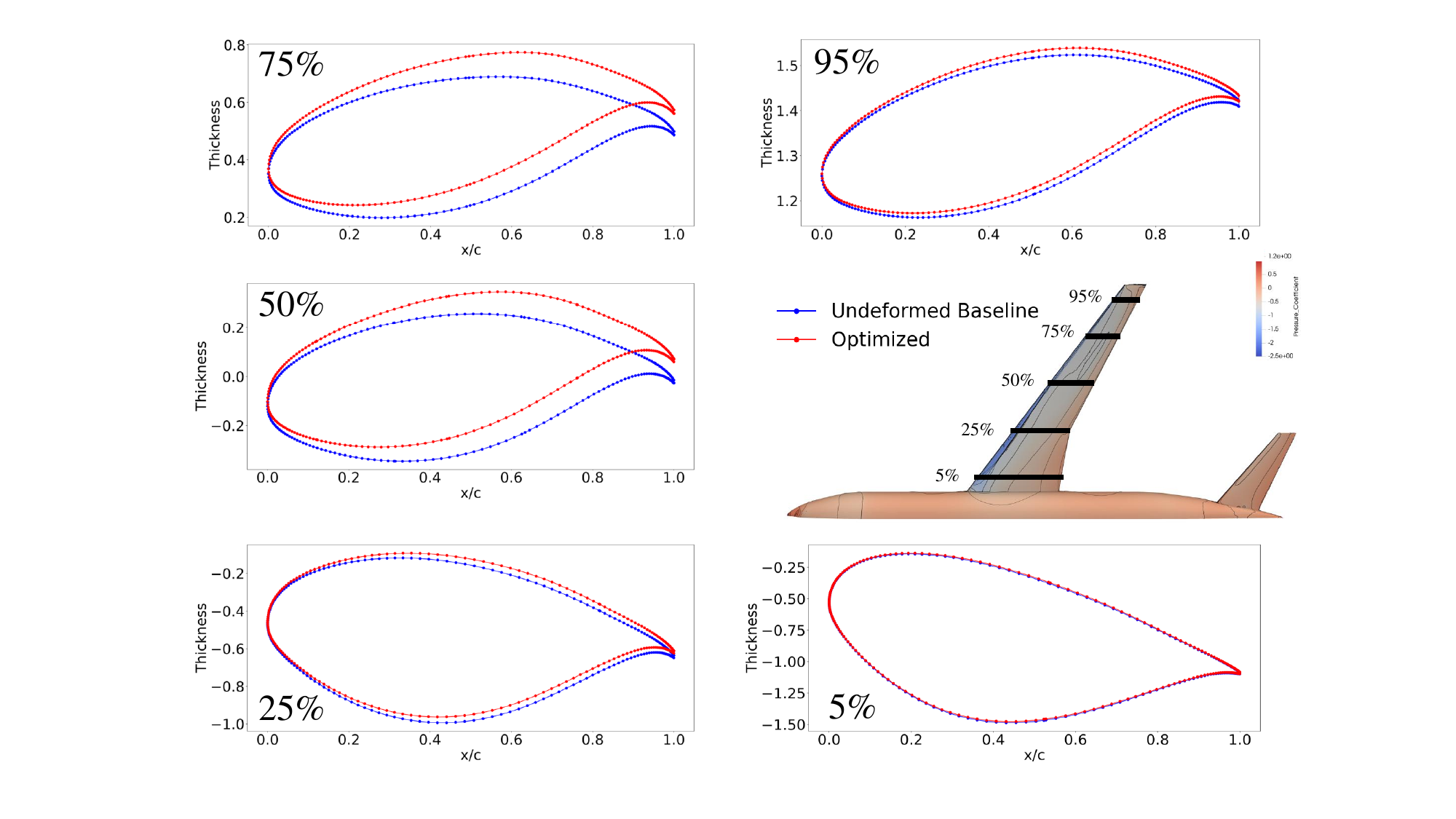}
    \caption{Comparison of airfoil shapes at various percentages of the span between an undeformed baseline wing and an optimized wing at a Pareto optimal point.}
    \label{fig:wing-shape-slices}
\end{figure}

Finally, we show the cross-sectional shapes for the Pareto optimal design versus the baseline design in \Cref{fig:wing-shape-slices}. Notice that at $5$ and $95$ \% wing span, the change in shape due to optimization is subtle, but the impact on the pressure distribution (\Cref{fig:cp-wing-sections}) is a lot more noticeable -- this demonstrates the sensitivity of the pressure coefficient to wing shape at the root and tip of the wing. In the remaining locations, the wing shape changes more noticeably, which includes a \mg{change in the overall wing twist}. In summary, multiobjective optimization identifies a design that offers a pragmatic compromise between the design performance at the two cruise points. Solving this problem efficiently and accurately was enabled by our proposed method \qpots.

\section{Conclusion}
\label{sec:conclusion}
Aircraft aerodynamic design optimization is better served by designing for the varying operating conditions at multiple sub-segments along the cruise segment, as opposed to considering a fixed operating condition. While existing work addresses this limitation by performing a multi-point optimization that considers a weighted combination of the objectives at all the sub-segments, such an approach depends on specifying appropriate weights and can miss out on identifying other, potentially more promising designs.
In contrast, this work establishes the necessity of solving a multiobjective optimization problem that simultaneously considers the objectives corresponding to all operating conditions. We showed that this is necessary to properly address the conflicts that typically exist between designs optimized individually for each sub-segment.

However, existing methods are prohibitive when applied to multiobjective optimization problems with expensive high-fidelity models. To address this, we propose a novel multiobjective Bayesian optimization method, called \qpots, that is inspired by Thompson sampling, involving Gaussian process surrogates. Specifically, our method adaptively selects a sequence of iterates that are Pareto optimal according to a posterior Gaussian process surrogate model. This way, our method hybridizes classical evolutionary algorithms---which are accurate but not sample efficient---with state-of-the-art Bayesian optimization---which is sample efficient but not accurate enough or scalable---to result in a more accurate and sample efficient approach. We demonstrate our approach on several synthetic experiments as well as the aerodynamic design optimization of the NASA common research model. In all the experiments, we show that our proposed approach either performs the best or is among the best performing.

\mg{\noindent {\bf Known limitations.} Despite the strengths of \qpots, we acknowledge a few limitations. First, the there is no automatic way to handle high-dimensional scaling which is important for aerodynamic design. However, high-dimensional scaling is a focused research direction in Bayesian optimization~\cite{rana2017high,zhang2019high,malu2021bayesian}. However, \qpots~is compatible with approaches that exist to address this, e.g., trust-region based multiobjective Bayesian optimization~\cite{daulton2022multi} -- this is a direction that we plan on pursuing in the future. Second, the posterior sampling complexity is well addressed via the Nystr\"{o}m approximation; however, a well-defined strategy to control the cost-accuracy tradeoff in the algorithm performance doesn't yet exist. We plan on addressing this as well in the future.
}

As a first step, this work considered the NASA CRM experiment with two objectives ($K=2$) and without constraints. Therefore, the first leg of future work will consider the inclusion of constraints and more sub-segments considered (higher $K$). We show in the Appendix that extensions to include constraints are straightforward in our approach. Another line of future work is to benchmark our method on the NASA CRM with batch acquisition $(q>1)$ to enable asynchronous parallel evaluations on high-performance computing environments.  This is necessary to minimize the wallclock times in the overall optimization involving expensive objectives and constraints. Finally, we also believe in exploring the utility of nonstationary surrogate models, e.g., deep Gaussian process models~\cite{rajaram2020deep,rajaram2021empirical,booth2024actively,booth2025contour}, within our framework. 
Our proposed framework is generic and replacing the surrogate model with another probabilistic surrogate model is straightforward.

\section*{Acknowledgements}
The authors of this work recognize the Penn State Institute for Computational and Data Sciences (RRID:SCR\_025154) for providing access to computational research infrastructure within the Roar Core Facility (RRID: SCR\_026424).

\clearpage
\section*{Appendix}

\subsection{Extension to constrained problems}
\label{sec:constraints}
 Our method seamlessly extends to constrained multiobjective problems, including both inequality and equality constraints, and the key step is explained as follows. We demonstrate this with the setting that we have a total of $C$ constraints $\{ c_1, \ldots, c_C\}$, and $\mcl{I} \cup \mcl{E} = [C]$. Further, let $Y_{k+i\in \mcl{I}}$ and $Y_{k+i\in \mcl{E}}$ denote the posterior GP of the $i$th inequality constraint and equality constraint functions, respectively. This way, we fit a total of $K+C$ posterior GPs $Y_1, \ldots, Y_K, Y_{K+1}, \ldots, Y_{K+C}$.
 For all $i \in [C]$, we define the constrained \qpots~ to choose points according to
\begin{equation}
\begin{split}
p_{\X^*}(\x) = \int \delta &\left(\x - \argmax_{\x \in \X}
\left[ \{Y_1(\x),\ldots,Y_K(\x)\} \times 
\mathds{1}_{\{\bigcap_{i} Y_{k+i \in \mcl{E}}(\x)= 0 \} }\times \mathds{1}_{\{\bigcap_i Y_{k+i \in \mcl{I}}(\x)\geq 0 \}}\right]\right) \\
& p(Y_1|\D_n)\ldots p(Y_K|\D_n)dY_1\ldots dY_K,
\end{split}
\end{equation}
where $\mathds{1}_{\bigcap_i Y_{k+i \in \mcl{I}}(\x)\geq 0}$ and $\mathds{1}_{\bigcap_i Y_{k+i \in \mcl{I}}(\x) = 0}$ are the indicator functions for joint feasibility with inequality and equality constraints, respectively, that take a value of $1$ when feasible and $0$ otherwise. This 
is equivalent to drawing samples from the GP posteriors, filtering them out by the indicator functions, and finding a Pareto optimal solution set per the posterior GP sample paths. 
\begin{equation}
    \begin{split}            
  X^* = \argmax_{\x \in \X} &\left[ \{Y_1(\x, \omega), \ldots, Y_K(\x, \omega)\}\times \mathds{1}_{\{\bigcap_{i} Y_{k+i \in \mcl{E}}(\x)= 0 \} }\times \mathds{1}_{\{\bigcap_i Y_{k+i \in \mcl{I}}(\x)\geq 0 \}}\right],
  \label{eqn:moo_gp_cons}
  \end{split}
\end{equation}
The rest of the \qpots~algorithm remains exactly the same.

\subsection{Proof that the multipoint optimum is a member of the Pareto frontier}
\label{sec:multip_v_multiob}
Let $\mathcal X$ be the domain of $f$ and $\bs{f}:\mathcal X\to\mathbb R^m$ be the vector of objectives to be minimized. Denote by
\[
F \;=\; \{\, f(\x) : \x\in\mathcal X\,\}\subset\mathbb R^m
\]
the attainable objective set. Let the definition for Pareto optimality presented in \Cref{sec:problem_statement} hold. Then, we proceed in two steps. First, we show that the weighted-sum minimizers (multi-point optima) are Pareto optimal. Then, we show that the converse is not necessarily always true.

\begin{proposition}
\label{thm:ws_is_pareto}
Let the weights in the multi-point optimum $w\in\mathbb R^m$ satisfy $w_i>0$ for every $i=1,\dots,m$. If
\[
\x^\star\in\arg\min_{\x\in\mathcal X} w^\top f(\x),
\]
then $\x^\star$ is Pareto optimal.
\end{proposition}

\begin{proof}
Assume, toward contradiction, that $\x^\star$ is not Pareto optimal. Then there exists $y\in\mathcal X$ with $f_i(\mbf{y})\le f_i(\x^\star)$ for all $i$ and $f_j(\mbf{y})<f_j(\x^\star)$ for some index $j$. Multiplying componentwise by the strictly positive weights $w_i$ and summing yields
\[
w^\top f(\mbf{y}) \;=\; \sum_{i=1}^m w_i f_i(\mbf{y})
\;<\; \sum_{i=1}^m w_i f_i(\x^\star) \;=\; w^\top f(\x^\star),
\]
which contradicts the minimality of $
x^\star$. Hence no such $\mbf{y}$ exists and $\x^\star$ must be Pareto optimal.
\end{proof}

\subsection*{Converse: when does every Pareto point arise from a weighted sum?}
The converse statement ("for every Pareto point there exists $w>0$ such that the point minimizes $w^\top f$") is \emph{not} true in general. A classical sufficient condition that guarantees the converse is convexity of the attainable objective set.

\begin{proposition}[Converse under convexity]
If the attainable objective set $F$ is convex, then for every Pareto optimal point $f^\star\in F$ there exists $w\in\mathbb R^m$ with $w_i\ge 0$, $w\ne 0$, such that $f^\star$ minimizes $w^\top y$ over $y\in F$. If additionally $f^\star$ is a \emph{strict} Pareto boundary point and appropriate regularity holds, one may choose $w_i>0$.
\end{proposition}

\begin{proof}[Sketch]
For a convex set $F\subset\mathbb R^m$ and a boundary point $f^\star\in\partial F$, the supporting hyperplane theorem guarantees the existence of a nonzero vector $w$ and scalar $\alpha$ with
\[
w^\top y \ge \alpha = w^\top f^\star\quad\forall y\in F.
\]
Thus $f^\star$ minimizes $w^\top y$ over $F$. Standard arguments translate this supporting hyperplane into a positive-weighted minimization if the Pareto point is not singular with respect to any objective (see multiobjective optimization texts for a full statement).
\end{proof}

Consequently, linear (weighted) scalarization is \emph{sufficient} but not \emph{necessary} to generate Pareto points unless $F$ is convex. Nonconvex portions of the Pareto front can therefore be missed by any linear scalarization.

\subsection{Impact of observation noise variance}
\label{sec:noise}

{We added a short robustness experiment in Section 5.3: for one synthetic test (Branin-Currin) we ran qPOTS with and without additive Gaussian observation noise (for a range of variance values). We report hypervolume convergence curves (mean $\pm$ std over 10 repeats) showing that qPOTS maintains stable convergence relative to baselines up to the tested noise levels.}
  \begin{figure}[t]
  \centering
  \includegraphics[width=1\linewidth]{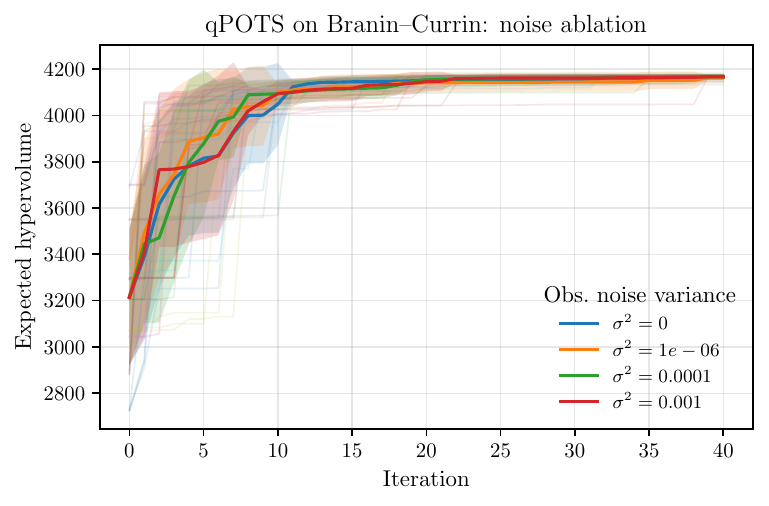}
  \caption{qPOTS on Branin--Currin under varying observation noise (\(\sigma^2\)).
  Curves show mean hypervolume across repetitions; shaded regions are \(\pm 1\) standard deviation.}
  \label{fig:branincurrin_noise_ablation}. 
\end{figure}

\clearpage

\bibliographystyle{plainnat}
\bibliography{references, multiobjective_aero_optimization}

\end{document}